\documentstyle{article}

\newtheorem{theorem}{Theorem}[section]
\newtheorem{corollary}[theorem]{Corollary}
\newtheorem{proposition}[theorem]{Proposition}
\newtheorem{lemma}[theorem]{Lemma}
\newtheorem{definition}[theorem]{Definition}

\begin{document}

\title{Q-Markov Random Probability Measures \\
and Their Posterior Distributions}

\author{R.M. Balan\thanks{Research done while the author was supported by a post-doctoral fellowship from the
 Natural Sciences and Engineering Research Council of Canada at Universit\'{e} de Sherbrooke.} \\
Department of Mathematics and Statistics \\
University of Ottawa \\
Ottawa, ON, K1N 6N5, Canada \\
{\em E-mail address}: rbala348@science.uottawa.ca }

\date{August 21, 2003}

\maketitle

\begin{abstract}
In this paper, we use the Markov property introduced in Balan and
Ivanoff (2002) for set-indexed processes and we prove that a
Markov prior distribution leads to a Markov posterior
distribution. In particular, by proving that a neutral to the
right prior distribution leads to a neutral to the right posterior
distribution, we extend a fundamental result of Doksum (1974) to
arbitrary sample spaces.
\end{abstract}



\section{Introduction}

Bayesian non-parametric statistics is a field that has been
introduced by Ferguson in 1973 and has become increasingly popular
among the theoretical statisticians in the past few decades. The
philosophy behind this field is to assume that the common
(unknown) distribution $P$ of a given sample
$\underline{X}=(X_{1}, \ldots, X_{n})$ is also governed by
randomness, and therefore can be regarded as a stochastic process
(indexed by sets). The best way for a Bayesian statistician to
guess the ``shape'' of the prior distribution $P$ is to identify
the posterior distribution of $P$ given $\underline{X}$ and to
prove that it satisfies the same properties as the prior.

Formalizing these ideas, we can say that a typical problem in
Bayesian nonparametric statistics is to identify a class $\Sigma$
of ``random distributions'' $P$ such that if $\underline{X}$ is a
sample of $n$ observations drawn according to $P$, then the
posterior distribution of $P$ given $\underline{X}$ remains in the
class $\Sigma$. The purpose of this paper is to introduce a new
class $\Sigma$ for which this property is preserved. This is the
class of ${\cal Q}$-Markov processes (or distributions), which
contains the extensively studied class of neutral to the right
processes.

There are two major contributions in the literature in this field.
The first one is Ferguson's (1973) fundamental paper where it is
shown that the posterior distribution of a Dirichlet process is
also Dirichlet. (By definition, a {\em Dirichlet process} with
parameter measure $\alpha$ has a Dirichlet finite dimensional
distribution with parameters $\alpha(A_{1}), \ldots,
\alpha(A_{k}), \alpha((\cup_{i=1}^{k}A_{i})^{c})$ over any
disjoint sets $A_{1}, \ldots,A_{k} \in {\cal B}$.) The second one
is Doksum's (1974) fundamental paper where it is proved that if
${\cal X}={\bf R}$, then the posterior distribution of a neutral
to the right process is also neutral to the right. (A random
probability distribution function $F:=(F_{t})_{t \in {\bf R}}$ is
{\em neutral to the right} if $F_{t_{1}},
(F_{t_{2}}-F_{t_{1}})/(1-F_{t_{1}}), \ldots,
(F_{t_{k}}-F_{t_{k-1}})/(1-F_{t_{k-1}})$ are independent $\forall
t_{1}< \ldots <t_{k}$, or equivalently, $Y_{t}:=- \ln(1-F_{t}),t
\in {\bf R}$ is a process with independent increments.) A quick
review of the literature to date (Ferguson, 1974; Ferguson and
Phadia, 1979; Dykstra and Laud, 1981; Hjort, 1990; Walker and
Muliere, 1997; Walker and Muliere, 1999) reveals that neutral to
the right processes have received considerably attention in the
past three decades, especially because of their appealing
representation using L\'{e}vy processes and because of their
applications in survival analysis, reliability theory, life
history data.

In the present paper we extend Doksum's result to the class of
${\cal Q}$-{\em Markov} processes introduced in Balan and Ivanoff
(2002), which are characterized by Markov-type finite dimensional
distributions. Unlike Doksum's paper (and unlike most of the
statistical papers generated by it) our results are valid for
arbitrary sample spaces ${\cal X}$, which can be endowed with a
certain topological structure (in particular for ${\cal X}={\bf
R}^{d}$). Our main result (Theorem \ref{main}) proves that if
$P:=(P_{A})_{A \in {\cal B}}$ is a set-Markov random probability
measure and $X_{1}, \ldots,X_{n}$ is a sample from $P$, then the
conditional distribution of $P$ given $X_{1}, \ldots,X_{n}$ is
also set-Markov. This result is new even in the case ${\cal
X}={\bf R}$, when the set-Markov property coincides with the
classical Markov property.

The paper is organized as follows:

In Section 2 we describe the structure that has to be imposed on
the sample space ${\cal X}$ (which will be assumed for the entire
paper); under this structure we identify the necessary ingredients
for the construction of set-Markov (respectively ${\cal
Q}$-Markov) random probability measure.

In Section 3 we introduce the Bayesian nonparametric framework and
we prove that a set-Markov prior distribution leads to a
set-Markov posterior distribution. The essence of all calculations
is an integral form of Bayes' formula.

In Section 4 we define neutral to the right processes and using
their ${\cal Q}$-Markov property we prove that a neutral to the
right prior distribution leads to a neutral to the right posterior
distribution.

The paper also includes two appendices: Appendix A contains two
elementary results which are used for the proof of Theorem
\ref{main}; Appendix B contains a Bayes property of a classical
Markov chain, which is interesting by itself and which has
motivated this paper.

\section{Q-Markov random probability measures}

Let $({\cal X}, {\cal B})$ be an arbitrary measurable space (the
sample space).

\begin{definition}
{\rm A collection $P:=(P_{A})_{A \in {\cal B}}$ of $[0,1]$-valued
random variables is called a {\em random probability measure} if

\noindent {\bf (i)} it is finitely additive in distribution, i.e.,
for every disjoint sets $(A_{j})_{j=1, \ldots,k}$ and for every $1
\leq i_{1}< \ldots <i_{m} \leq k$, the distribution of
$(P_{\cup_{j=1}^{i_{1}}A_{j}},
\ldots,P_{\cup_{j=i_{m}}^{k}A_{j}})$ coincides with the
distribution of $(\sum_{j=1}^{i_{1}}P_{A_{j}}, \ldots,
\sum_{j=i_{m}}^{k}P_{A_{j}})$;

\noindent {\bf (ii)} $P_{\cal X}=1$ a.s.; and

\noindent {\bf (iii)} it is countably additive in distribution,
i.e., for every decreasing sequence $(A_{n})_{n} \subseteq {\cal
B}$ with $\cap_{n}A_{n}=\emptyset$ we have $\lim_{n}P_{A_{n}}=0$
a.s.}
\end{definition}

Note that the almost sure convergence of {\bf (iii)} (in the above
definition) is equivalent to the convergence in distribution and
the convergence in mean.

In order to construct a random probability measure $P$ on ${\cal
B}$ it is enough to specify its finite dimensional distributions
$\mu_{A_{1} \ldots A_{k}}$ over all un-ordered collections
$\{A_{1}, \ldots, A_{k}\}$ of disjoint sets in ${\cal B}$. Some
conditions need to be imposed.

\vspace{2mm}

{\em Condition C1.} If $\{A_{1}, \ldots,A_{k}\}$ is an un-ordered
collection of disjoint sets and we let
$A'_{l}:=\cup_{j=i_{l-1}+1}^{i_{l}}A_{j};l=1, \ldots,m$ for $1
\leq i_{1}< \ldots <i_{m} \leq k$, then $\mu_{A'_{1} \ldots
A'_{m}}=\mu_{A_{1} \ldots A_{k}} \circ \alpha^{-1}$, where
$\alpha(x_{1}, \ldots, x_{k})=(\sum_{j=1}^{i_{1}}x_{j}, \ldots,
\sum_{j=i_{m-1}+1}^{i_{m}}x_{j})$.

\vspace{2mm}

{\em Condition C2.} For every $(A_{n})_{n} \subseteq {\cal B}$
with $A_{n+1} \subseteq A_{n}, \forall n$ and
$\cap_{n}A_{n}=\emptyset$, we have
$\lim_{n}\mu_{A_{n}}=\delta_{0}$.

\vspace{2mm}

 In this paper we will assume that the sample space
${\cal X}$ has an additional underlying structure which we begin
now to explain.

Let ${\cal X}$ be a (Hausdorff) topological space and ${\cal B}$
its Borel $\sigma$-field. We will assume that there exists a
collection ${\cal A}$ of closed subsets of ${\cal X}$ which
generates ${\cal B}$ (i.e. ${\cal B}=\sigma({\cal A})$)  and which
has the following properties:
\begin{enumerate}

\item $\emptyset,{\cal X} \in {\cal A}$; \item ${\cal A}$ is a
semilattice i.e., ${\cal A}$ is closed under arbitrary
intersections;

\item $\forall A,B \in {\cal A}; A,B \not = \emptyset \Rightarrow
A \cap B \not = \emptyset$;

\item There exists a sequence $({\cal A}_{n})_{n}$ of finite
sub-semilattices of ${\cal A}$ such that $\forall A \in {\cal A}$,
there exist $A_{n} \in {\cal A}_{n}(u), \forall n$ with
$A=\cap_{n}A_{n}$ and $A \subseteq A_{n}^{0}, \forall n$. (Here
${\cal A}_{n}(u)$ denotes the class of all finite unions of sets
in ${\cal A}_{n}$.)
\end{enumerate}

More details about this type of structure can be found in Ivanoff
and Merzbach (2000), where ${\cal A}$ is called an {\em indexing
collection}. By properties 2 and 3, the collection ${\cal A}$ has
the finite intersection property, and hence its minimal set
$\emptyset':= \cap_{A \in {\cal A} \verb2\2 \{\emptyset\}}A$ is
non-empty.

The typical example of a sample space ${\cal X}$ which can be
endowed with an indexing collection is ${\bf R}^{d}$; in this case
${\cal A}=\{[0,z];z \in {\bf R}^{d}\} \cup \{\emptyset, {\bf
R}^{d}\}$ and the approximation sets $A_{n}$ have vertices with
dyadic coordinates.

We denote with ${\cal A}(u)$ the class of all finite unions of
sets in ${\cal A}$, with ${\cal C}$ the semialgebra of the sets
$C=A \verb2\2 B$ with $A \in {\cal A},B \in {\cal A}(u)$
 and with ${\cal C}(u)$ the algebra of sets generated by ${\cal C}$. Note that ${\cal B}=\sigma({\cal C}(u))$.

We introduce now the definition of the ${\cal Q}$-Markov property.
This definition has been originally considered in Balan and
Ivanoff (2002) for finitely additive real-valued processes indexed
by the algebra ${\cal C}(u)$. In this paper, we will restrict our
attention to random probability measures.

\begin{definition}
\label{definition-Q} {\bf (a)} For each $B_{1},B_{2} \in {\cal
A}(u)$ with $B_{1} \subseteq B_{2}$, let $Q_{B_{1}B_{2}}$ be a
transition probability on $[0,1]$. The family ${\cal
Q}:=(Q_{B_{1}B_{2}})_{B_{1} \subseteq B_{2}}$ is called a {\bf
transition system} if $\forall B_{1} \subseteq B_{2} \subseteq
B_{3}$ in ${\cal A}(u), \forall z_{1} \in [0,1], \forall
\Gamma_{3} \in {\cal B}([0,1])$
$$Q_{B_{1}B_{3}}(z_{1}; \Gamma_{3})=\int_{[0,1]}Q_{B_{2}B_{3}}(z_{2};\Gamma_{3})Q_{B_{1}B_{2}}(z_{1};dz_{2})$$

\noindent {\bf (b)} Given a transition system ${\cal
Q}:=(Q_{B_{1}B_{2}})_{B_{1} \subseteq B_{2}}$, a random
probability measure $P:=(P_{A})_{A \in {\cal B}}$, defined on a
probability space $(\Omega, {\cal F}, {\cal P})$, is called {\bf
${\cal Q}$-Markov} if $\forall B_{1} \subseteq B_{2}$ in ${\cal
A}(u)$, $\forall \Gamma_{2} \in {\cal B}([0,1])$
$${\cal P}[P_{B_{2}} \in \Gamma_{2}|{\cal F}_{B_{1}}]=Q_{B_{1}B_{2}}(P_{B_{1}};
\Gamma_{2}) \ \ {\rm a.s.}$$ where ${\cal
F}_{B_{1}}:=\sigma(\{P_{A}; A \in {\cal A}, A \subseteq B_{1}\})$.
\end{definition}

A ${\cal Q}$-Markov random probability measure can be constructed
using the following additional consistency condition.

\vspace{3mm}

{\em Condition C3.} If $(Y_{1}, \ldots, Y_{k})$ is a vector with
distribution $\mu_{C_{1} \ldots C_{k}}$ where
$C_{1}=B_{1};C_{i}=B_{i} \verb2\2 B_{i-1};i=2, \ldots,k$ and
$B_{1} \subseteq \ldots \subseteq B_{k}$ are sets in ${\cal
A}(u)$, then for every $i=2, \ldots,k$, the distribution of
$Y_{i}$ given $Y_{1}=y_{1}, \ldots, Y_{i-1}=y_{i-1}$ depends only
on $y:=\sum_{j=1}^{i-1}y_{j}$ and is equal to $Q_{B_{i-1}B_{i}}(y;
y+ \cdot)$.

\vspace{3mm}

\noindent The next result follows immediately by Kolmogorov's
extension theorem.

\begin{theorem}
\label{constr-set-Markov} Let ${\cal Q}:=(Q_{B_{1}B_{2}})_{B_{1}
\subseteq B_{2}}$ be a transition system. For each un-ordered
collection $\{A_{1}, \ldots,A_{k}\}$ of disjoint sets in ${\cal
B}$ let $\mu_{A_{1} \ldots A_{k}}$ be a probability measure on
$([0,1]^{k}, {\cal B}([0,1])^{k})$
 such that C1-C3 hold; let $\mu_{\emptyset}=\delta_{0}, \mu_{\cal X}=\delta_{1}$.
Then there exists a probability measure ${\cal P}^{1}$ on
$([0,1]^{\cal B}, {\cal B}([0,1])^{\cal B})$ under which the
coordinate-variable process $P:=(P_{A})_{A \in {\cal B}}$ is a
${\cal Q}$-Markov random probability measure whose finite
dimensional distributions are the measures $\mu_{A_{1} \ldots
A_{k}}$.
\end{theorem}

{\bf Examples}:
\begin{enumerate}

\item Let $P$ be the Dirichlet process with parameter measure
$\alpha$. For any disjoint sets $A_{1}, \ldots, A_{k}$ in ${\cal
B}$, $(P_{A_{1}}, \ldots, P_{A_{k}})$ has a {\em Dirichlet}
distribution with parameters $\alpha(A_{1}), \ldots,\alpha(A_{k}),
\alpha((\cup_{i=1}^{k}A_{i})^{c})$. The ratio
$P_{A_{i}}/(1-\sum_{j=1}^{i-1}P_{A_{j}})$ is independent of
$P_{A_{1}}, \ldots, P_{A_{i-1}}$ and has a Beta distribution with
parameters $\alpha(A_{i}), \alpha((\cup_{j=1}^{i}A_{j})^{c})$;
hence the distribution of $P_{A_{i}}$ given $P_{A_{1}}, \ldots,
P_{A_{i-1}}$ depends only on $\sum_{j=1}^{i-1}P_{A_{j}}$. The
process $P$ is ${\cal Q}$-Markov with $Q_{B_{1}B_{2}}(z_{1};
\Gamma_{2})$ equal to the value at $(\Gamma_{2}-z_{1})/(1-z_{1})$
of the Beta distribution with parameters $\alpha(B_{2} \verb2\2
B_{1}),\alpha(B_{2}^{c})$.

\item Let $P:=(1/N) \sum_{j=1}^{N}\delta_{Z_{j}}$ be the empirical
measure of a sample $Z_{1}, \ldots ,Z_{N}$ from a non-random
distribution $P_{0}$ on ${\cal X}$. For any disjoint sets $A_{1},
\ldots, A_{k}$ in ${\cal B}$, $(NP_{A_{1}}, \ldots, NP_{A_{k}})$
has a {\em multinomial} distribution with $N$ trials and
$P_{0}(A_{1}), \ldots, P_{0}(A_{k})$ probabilities of success;
hence the distribution of $NP_{A_{i}}$ given $NP_{A_{1}}, \ldots,
NP_{A_{i-1}}$ depends only on $\sum_{j=1}^{i-1}P_{A_{j}}$ (it is a
binomial distribution with $N(1- \sum_{j=1}^{i-1}P_{A_{j}})$
trials and $P_{0}(A_{i})/(1-\sum_{j=1}^{i-1}P_{0}(A_{j}))$
probability of success). The process $P$ is ${\cal Q}$-Markov with
$$Q_{B_{1}B_{2}}\left(\frac{m_{1}}{N}; \left\{ \frac{m_{2}}{N} \right\}\right)=
\left( \begin{array}{c}N-m_{1} \\ m_{2}-m_{1} \end{array} \right)
 \frac{P_{0}(C)^{m_{2}-m_{1}}P_{0}(B_{2}^{c})^{N-m_{2}}}{P_{0}(B_{1}^{c})^{N-m_{1}}}$$
where $\left( \begin{array}{c} a \\ b \end{array}
\right)=a!/b!(a-b)!$ is the binomial coefficient and $C=B_{2}
\verb2\2 B_{1}$.

\item Let $P:=(1/N) \sum_{j=1}^{N}\delta_{W_{j}}$ be the empirical
measure of a sample $W_{1}, \ldots ,W_{N}$ from a Dirichlet
process with parameter measure $\alpha$. For any disjoint sets
$A_{1}, \ldots, A_{k}$ in ${\cal B}$, $(NP_{A_{1}}, \ldots,
NP_{A_{k}})$ has a {\em P\'{o}lya} distribution with $N$ trials
and parameters $\alpha(A_{1}), \ldots,\alpha(A_{k}),
\alpha((\cup_{i=1}^{k}A_{i})^{c})$; hence the distribution of
$NP_{A_{i}}$ given $NP_{A_{1}}, \ldots, NP_{A_{i-1}}$ depends only
on $\sum_{j=1}^{i-1}P_{A_{j}}$ (it is a P\'{o}lya distribution
with $N(1- \sum_{j=1}^{i-1}P_{A_{j}})$ trials and parameters
$\alpha(A_{i}), \alpha((\cup_{j=1}^{i}A_{j})^{c})$). The process
$P$ is ${\cal Q}$-Markov with
$$Q_{B_{1}B_{2}}\left(\frac{m_{1}}{N}; \left\{ \frac{m_{2}}{N} \right\}\right)=
\left( \begin{array}{c}N-m_{1} \\ m_{2}-m_{1} \end{array} \right)
 \frac{\alpha(C)^{[m_{2}-m_{1}]}\alpha(B_{2}^{c})^{[N-m_{2}]}}{\alpha(B_{1}^{c})^{[N-m_{1}]}}$$
where $\alpha^{[x]}=\alpha (\alpha+1) \ldots (\alpha+x-1)$  and
$C=B_{2} \verb2\2 B_{1}$.

\end{enumerate}

\section{The posterior distribution of a Q-Markov random probability measure}

We begin to introduce the Bayesian nonparametric framework.

Let $P:=(P_{A})_{A \in {\cal B}}$ be a ${\cal Q}$-Markov random
probability measure defined on a probability space $(\Omega, {\cal
F}, {\cal P})$ and $X_{i}: \Omega \rightarrow {\cal X}, i=1,
\ldots,n$ some ${\cal F}/{\cal B}$-measurable functions such that
$\forall A_{1}, \ldots,A_{n} \in {\cal B}$
$${\cal P}[X_{1} \in A_{1}, \ldots,X_{n} \in A_{n}|P]=\prod_{i=1}^{n}P_{A_{i}} \ \ {\rm a.s.}$$

We say that $\underline{X}:=(X_{1}, \ldots,X_{n})$ is a {\em
sample from $P$}. The distribution of $P$ is called {\em prior},
while the distribution of $P$ given $\underline{X}$ is called {\em
posterior}. Note that $(P_{A})_{A \in {\cal B}}$ and $X_{1},
\ldots, X_{n}$ can be constructed as coordinate-variables on the
space $([0,1]^{\cal B} \times {\cal X}^{n}, {\cal B}([0,1])^{\cal
B} \times {\cal B}^{n})$ under the probability measure ${\cal P}$
defined by
$${\cal P}(D \times \prod_{i=1}^{n}A_{i}):=\int_{D} \prod_{i=1}^{n}\omega_{A_{i}} \ {\cal P}^{1}(d\omega), \ \ D \in {\cal B}([0,1])^{\cal B}, A_{i} \in {\cal B}$$
where ${\cal P}^{1}$ is the probability measure given by Theorem
\ref{constr-set-Markov}.

The goal of this section is to prove that the posterior
distribution of $P$ given $\underline{X}={\underline x}$ is ${\cal
Q}^{({\underline x})}$-Markov (for some ``posterior'' transition
system ${\cal Q}^{({\underline x})}$).

\vspace{3mm}

Let $\alpha_{n}$ be the law of $\underline{X}$ under ${\cal P}$
and $\mu_{A_{1}, \ldots, A_{k}}$ be the law of $(P_{A_{1}},
\ldots, P_{A_{k}})$ under ${\cal P}$, for every $A_{1}, \ldots,
A_{k} \in {\cal B}$. Note that
$\alpha_{n}(\prod_{i=1}^{n}A_{i})={\cal
E}[\prod_{i=1}^{n}P_{A_{i}}]$, where ${\cal E}$ denotes the
expectation with respect to ${\cal P}$.

 For each set $B_{1} \in
{\cal A}(u)$, let $\nu_{B_{1}}$ be the law of $(X_{1}, \ldots,
X_{n},P_{B_{1}})$ under ${\cal P}$. Note that
$\nu_{B_{1}}(\prod_{i=1}^{n}A_{i} \times \Gamma_{1})={\cal
E}[\prod_{i=1}^{n}P_{A_{i}} \cdot I_{\Gamma_{1}}(P_{B_{1}})]$ and
\begin{equation}
\label{disint-nu-B1-n} \nu_{B_{1}}(\tilde{A} \times \Gamma_{1}) =
\int_{\tilde{A}} \mu_{B_{1}}^{(\underline{x})} (\Gamma_{1})
\alpha_{n}(d \underline{x}) = \int_{\Gamma_{1}}
\tilde{Q}_{B_{1}}(z_{1};\tilde{A}) \mu_{B_{1}}(dz_{1})
\end{equation}
where $\mu_{B_{1}}^{(\underline{x})}(\Gamma_{1}):={\cal
P}[P_{B_{1}} \in \Gamma_{1}|\underline{X}=\underline{x}]$ and
$\tilde{Q}_{B_{1}}(z_{1};\tilde{A}):={\cal P}[\underline{X} \in
\tilde{A}|P_{B_{1}}=z_{1}]$.

\vspace{2mm}

For each sets $B_{1},B_{2} \in {\cal A}(u);B_{1} \subseteq B_{2}$,
let $\nu_{B_{1}B_{2}}$ be the law of $(X_{1}, \ldots, X_{n},
 \linebreak P_{B_{1}}, P_{B_{2}})$ under ${\cal P}$. Note that
$\nu_{B_{1}B_{2}}(\prod_{i=1}^{n}A_{i} \times \Gamma_{1} \times
\Gamma_{2}) = {\cal E}[\prod_{i=1}^{n}P_{A_{i}} \cdot
I_{\Gamma_{1}}(P_{B_{1}}) \linebreak I_{\Gamma_{2}}(P_{B_{2}})]$
and
\begin{eqnarray}
\label{disint-nu-B1-B2-n} \nu_{B_{1}B_{2}}(\tilde{A} \times
\Gamma_{1} \times \Gamma_{2}) & = & \int_{\tilde{A}}
\int_{\Gamma_{1}}
Q_{B_{1}B_{2}}^{(\underline{x})}(z_{1};\Gamma_{2})
\mu_{B_{1}}^{(\underline{x})}(dz_{1})
\alpha_{n}(d \underline{x}) \\
 & = & \int_{\Gamma_{1} \times \Gamma_{2}} \tilde{Q}_{B_{1}B_{2}}(z_{1},z_{2}; \tilde{A}) \mu_{B_{1}B_{2}}(dz_{1} \times dz_{2})
\end{eqnarray}
where
\begin{equation}
\label{definition-Q-x}
 Q_{B_{1}B_{2}}^{(\underline{x})}(z_{1}; \Gamma_{2}):={\cal
P}[P_{B_{2}} \in \Gamma_{2}|\underline{X}=
\underline{x},P_{B_{1}}=z_{1}] \end{equation}
 and
$\tilde{Q}_{B_{1}B_{2}}(z_{1},z_{2}; \tilde{A}):= {\cal
P}[\underline{X} \in \tilde{A}|P_{B_{1}}=z_{1},P_{B_{2}}=z_{2}]$.
(For the first equality we used the first integral in the
decomposition (\ref{disint-nu-B1-n}) of $\nu_{B_{1}}$).

\noindent Using the second integral in the decomposition
(\ref{disint-nu-B1-n}) of $\nu_{B_{1}}$ and the ${\cal Q}$-Markov
property for representing $\mu_{B_{1}B_{2}}$ we get: (for
$\mu_{B_{1}}$-almost all $z_{1}$)
\begin{equation}
\label{central-Bayes}
\int_{\tilde{A}}Q_{B_{1}B_{2}}^{(\underline{x})}(z_{1};\Gamma_{2})
\tilde{Q}_{B_{1}}(z_{1};d \underline{x})= \int_{\Gamma_{2}}
\tilde{Q}_{B_{1}B_{2}}(z_{1},z_{2};
\tilde{A})Q_{B_{1}B_{2}}(z_{1};dz_{2}).
\end{equation}

This very important equation is the key for determining the
posterior transition probabilities
$Q_{B_{1}B_{2}}^{(\underline{x})}$ from the prior transition
probabilities $Q_{B_{1}B_{2}}$, providing that
$\tilde{Q}_{B_{1}}(z_{1}; \prod_{i=1}^{n}A_{i})={\cal
E}[\prod_{i=1}^{n}P_{A_{i}}| P_{B_{1}}=z_{1}]$ and
$\tilde{Q}_{B_{1}B_{2}}(z_{1},z_{2}; \prod_{i=1}^{n}A_{i}) = {\cal
E}[\prod_{i=1}^{n}P_{A_{i}}| P_{B_{1}}=z_{1},P_{B_{2}}=z_{2}]$ are
easily computable.

We note that each $Q_{B_{1}B_{2}}^{(x)}(z_{1}; \cdot)$ is
well-defined only for $\nu_{B_{1}}$-almost all $(\underline{x},
z_{1})$. Moreover, as we will see in the proof of Theorem
\ref{main} and it was correctly pointed out by an anonymous
referee, ${\cal Q}^{(\underline{x})}$ may not be a genuine
transition system as introduced by Definition
\ref{definition-Q}.(a). To avoid any confusion we introduce the
following terminology.

\begin{definition}
The family ${\cal Q}^{(x)}:=(Q_{B_{1}B_{2}}^{(x)})_{B_{1}
\subseteq B_{2}}$ defined by (\ref{definition-Q-x}) is called a
{\bf posterior transition system} (corresponding to $P$ and
$\underline{X}$) if $\forall B_{1} \subseteq B_{2} \subseteq
B_{3}$ in ${\cal A}(u)$, $\forall \Gamma_{3} \in {\cal B}([0,1])$
and for $\nu_{B_{1}}$-almost all $(\underline{x},z_{1})$
$$Q_{B_{1}B_{3}}^{(\underline{x})}(z_{1}; \Gamma_{3})= \int_{[0,1]}
Q_{B_{2}B_{3}}^{(\underline{x})}(z_{2};
\Gamma_{3})Q_{B_{1}B_{2}}^{(\underline{x})}(z_{1};dz_{2})$$ In
this case, we will say that the conditional distribution of $P$
given $\underline{X}=\underline{x}$ is {\bf ${\cal
Q}^{(\underline{x})}$-Markov} if $\forall B_{1} \subseteq B_{2}$
in ${\cal A}(u)$, $\forall \Gamma_{2} \in {\cal B}([0,1]$
$${\cal P}[P_{B_{2}} \in \Gamma_{2}|{\cal F}_{B_{1}},\underline{X}]=Q_{B_{1}B_{2}}^{(\underline{X})}(P_{B_{1}};
\Gamma_{2}) \ \ {\rm a.s.}$$
\end{definition}

\vspace{3mm}

We proceed now to the proof of the main theorem. Two preliminary
lemmas are needed.

Let $B_{1} \subseteq B_{2}$ be some arbitrary sets in ${\cal
A}(u)$, $C:=B_{2} \verb2\2 B_{1}$ and $0 \leq l \leq r \leq n$.
The next lemma shows us what happens intuitively with the
probability that the first $l$ observations fall in $B_{1}$, the
next $r-l$ observations fall in $C$ and the remaining $n-r$
observations fall in $B_{2}^{c}$, given $P_{B_{1}}$ and
$P_{B_{2}}$.

\begin{lemma}
\label{Qtilde-disintegration} For each $B_{1} \subseteq B_{2}$ in
${\cal A}(u)$ and $A_{1}, \ldots, A_{n} \in {\cal B}$, let
\begin{equation}
\label{tilde-A} \tilde{A}:=\prod_{i=1}^{l}(A_{i} \cap B_{1})
\times \prod_{i=l+1}^{r}(A_{i} \cap C) \times
\prod_{i=r+1}^{n}(A_{i} \cap B_{2}^{c})
\end{equation}
where $C:=B_{2} \verb2\2 B_{1}$ and $0 \leq l \leq r \leq n$. Let
$\tilde{A}_{1}:=\prod_{i=1}^{l}(A_{i} \cap B_{1}) \times {\cal
X}^{n-l}$, $\tilde{A}_{2}:= \prod_{i=l+1}^{r}(A_{i} \cap C) \times
{\cal X}^{n-r+l}$, $\tilde{A}_{3}:= \prod_{i=r+1}^{n}(A_{i} \cap
B_{2}^{c}) \times {\cal X}^{r}$, $\tilde A_{23}:=\tilde A_{2} \cap
\tilde A_{3}$.

(a) For $\mu_{B_{1}}$-almost all $z_{1}$,
$\tilde{Q}_{B_{1}}(z_{1};
\tilde{A})=\tilde{Q}_{B_{1}}(z_{1};\tilde{A}_{1}) \cdot
\tilde{Q}_{B_{1}}(z_{1}; \tilde{A}_{23})$.

(b) For $\mu_{B_{1}B_{2}}$-almost all $(z_{1},z_{2})$,
$$\tilde{Q}_{B_{1}B_{2}}(z_{1},z_{2}; \tilde{A})=\tilde{Q}_{B_{1}}(z_{1};\tilde{A}_{1}) \cdot
\tilde{Q}_{B_{1}B_{2}}(z_{1},z_{2}; \tilde{A}_{2}) \cdot
\tilde{Q}_{B_{2}}(z_{2}; \tilde{A}_{3}).$$
\end{lemma}

\noindent {\bf Proof}: We will prove only (b) since part (a)
follows by a similar argument. Note that the sets $\tilde A$ form
a $\pi$-system generating the $\sigma$-field ${\cal B}^{n}$ on
 $B_{1}^{l} \times C^{r-l} \times (B_{2}^{c})^{n-r}$.

Since $\sigma({\cal A})={\cal B}$ and ${\cal A}$ is a
$\pi$-system, using a Dynkin system argument, it is enough to
consider the case $A_{1}, \ldots,A_{n} \in {\cal A}$. Note that
$${\cal E}[\prod_{i=r+1}^{n}P_{A_{i} \cap B_{2}^{c}}|{\cal F}_{B_{2}}] =
{\cal E}[\prod_{i=r+1}^{n}P_{A_{i} \cap B_{2}^{c}}|P_{B_{2}}]=
\tilde{Q}_{B_{2}}(P_{B_{2}};\tilde{A}_{3}).$$

\noindent By double conditionning with respect to ${\cal
F}_{B_{2}}$, we have
$$\tilde{Q}_{B_{1}B_{2}}(z_{1},z_{2}; \tilde{A})={\cal E}[\prod_{i=1}^{l}P_{A_{i} \cap B_{1}} \prod_{i=l+1}^{r}P_{A_{i} \cap C} \prod_{i=r+1}^{n}P_{A_{i} \cap B_{2}^{c}} \ | \ P_{B_{1}}=z_{1},P_{B_{2}}=z_{2}]=$$
$$\tilde{Q}_{B_{2}}(z_{2};\tilde{A}_{3}) \cdot {\cal E}[\prod_{i=1}^{l}P_{A_{i} \cap B_{1}} \cdot \prod_{i=l+1}^{r}P_{A_{i} \cap C} \ | \ P_{B_{1}}=z_{1},P_{B_{2}}=z_{2}].$$

\noindent For the second term we have
$${\cal E}[\prod_{i=1}^{l}P_{A_{i} \cap B_{1}} \prod_{i=l+1}^{r}P_{A_{i} \cap C} \ | \ P_{B_{1}},P_{B_{2}}]=$$
$${\cal E}[\prod_{i=1}^{l}P_{A_{i} \cap B_{1}} {\cal E}[\prod_{i=l+1}^{r}P_{A_{i} \cap C}|(P_{A_{i} \cap B_{1}})_{i \leq l},P_{B_{1}},P_{B_{2}}] \ | \ P_{B_{1}},P_{B_{2}}].$$

\noindent Since $P_{A_{i} \cap C}=P_{B_{1} \cup (A_{i} \cap
B_{2})}-P_{B_{1}}$, using Lemma \ref{lemmaA1} (Appendix A)
$${\cal E}[\prod_{i=l+1}^{r}P_{A_{i} \cap C}|(P_{A_{i} \cap B_{1}})_{i \leq l},P_{B_{1}},P_{B_{2}}]=\tilde{Q}_{B_{1}B_{2}}(P_{B_{1}},P_{B_{2}}; \tilde{A}_{2}).$$

\noindent (In order to use Lemma \ref{lemmaA1}, we need $A_{l+1}
\subseteq A_{l+2} \subseteq \ldots \subseteq A_{r}$. Note that
this is not a restriction since if we can consider the minimal
semilattice $\{A'_{1},\ldots,A'_{m}\}$ determined by the sets
$A_{l+1}, \ldots,A_{r}$, which is ordered such that $A'_{j} \not
\subseteq \cup_{l \not = j}A'_{l} \forall j$, and we let
$B'_{j}=\cup_{s=1}^{j}A'_{s}$ and $C'_{j}=B'_{j} \verb2\2
B'_{j-1}$, then each $A_{i}= \dot \cup_{j \in J_{i}}C'_{j}$ for
some $J_{i} \subseteq \{1, \ldots,m\}$. We have $A_{i} \cap C=\dot
\cup_{j \in J_{i}}[(B'_{j} \cap C) \verb2\2 (B'_{j-1} \cap C)]$
and $\prod_{i=l+1}^{r}P_{A_{i} \cap C}=h(P_{B'_{1} \cap C},
\ldots, P_{B'_{m} \cap C})$ for some function $h$.)

\noindent Finally, since ${\cal F}_{B_{1}}$ is conditionally
independent of $P_{B_{2}}$ given $P_{B_{1}}$ and $P_{A_{i} \cap
B_{1}}, i \leq l$ are ${\cal F}_{B_{1}}$-measurable, we have
${\cal E}[\prod_{i=1}^{l}P_{A_{i} \cap B_{1}} \ | \
P_{B_{1}},P_{B_{2}}] = {\cal E}[\prod_{i=1}^{l}P_{A_{i} \cap
B_{1}} \ | \ P_{B_{1}}] \linebreak
=\tilde{Q}_{B_{1}}(P_{B_{1}};\tilde{A}_{1})$, which concludes the
proof. $\Box$

\vspace{3mm}

\noindent {\em Note}: Let $\tilde A_{12}:=\tilde A_{1} \cap \tilde
A_{2}$. By a similar argument one can show that
\begin{equation}
\label{Q-tilda-B1-B2'} \tilde{Q}_{B_{1}B_{2}}(z_{1},z_{2};
\tilde{A}_{12})= \tilde{Q}_{B_{1}}(z_{1};\tilde{A}_{1}) \cdot
\tilde{Q}_{B_{1}B_{2}}(z_{1},z_{2}; \tilde{A}_{2})
\end{equation}
\begin{equation}
\label{Q-tilda-B1-B2''} \tilde{Q}_{B_{1}B_{2}}(z_{1},z_{2};
\tilde{A}_{23})= \tilde{Q}_{B_{1}B_{2}}(z_{1},z_{2};\tilde{A}_{2})
\cdot \tilde{Q}_{B_{2}}(z_{2};\tilde{A}_{3})
\end{equation}

The next lemma tells us that if $B_{1} \subseteq B_{2}$ are
``nicely-shaped'' regions and we want to predict the value of
$P_{B_{2}}$ given the value of $P_{B_{1}}$ and a sample
$\underline{X}$ from $P$, then we can forget all about those
values $X_{i}$ which fall inside the region $B_{1}$. The reason
for this phenomenon is the very essence of the Markov property
given by Definition \ref{definition-Q}.(b), which says that for
predicting the value of $P_{B_{2}}$ it suffices to know the value
of $P_{B_{1}}$, i.e. all the information about the values of $P$
inside the region $B_{1}$ can be discarded.

\begin{lemma}
For every $B_{1},B_{2} \in {\cal A}(u)$ with $B_{1} \subseteq
B_{2}$, for every $\Gamma_{2} \in {\cal B}([0,1])$ and for
$\nu_{B_{1}}$-almost all $(\underline{x},z_{1})$,
$Q_{B_{1}B_{2}}^{(\underline{x})}(z_{1}; \Gamma_{2})$ does not
depend on those $x_{i}$'s that fall in $B_{1}$; in particular, for
$\nu_{B_{1}}$-almost all $(\underline{x},z_{1})$ in $B_{1}^{n}
\times [0,1]$, $Q_{B_{1}B_{2}}^{(\underline{x})}(z_{1};
\Gamma_{2})=Q_{B_{1}B_{2}}(z_{1}; \Gamma_{2})$.
\end{lemma}

\noindent {\bf Proof}: Let $A_{1}, \ldots, A_{n} \in {\cal B}$ and
$\tilde{A}$ defined by (\ref{tilde-A}). Using
(\ref{central-Bayes}) and Lemma \ref{Qtilde-disintegration},(b)
combined with (\ref{Q-tilda-B1-B2''}) we have
$$\int_{\tilde{A}} Q_{B_{1}B_{2}}^{(\underline{x})}(z_{1};\Gamma_{2}) \tilde{Q}_{B_{1}}(z_{1};d \underline{x}) =\int_{\Gamma_{2}} \tilde{Q}_{B_{1}B_{2}}(z_{1},z_{2};\tilde{A})Q_{B_{1}B_{2}}(z_{1};dz_{2})=$$
$$\tilde{Q}_{B_{1}}(z_{1};\tilde{A}_{1}) \int_{\Gamma_{2}} \tilde{Q}_{B_{1}B_{2}}(z_{1},z_{2};\tilde{A}_{23})Q_{B_{1}B_{2}}(z_{1};dz_{2})=$$
$$\tilde{Q}_{B_{1}}(z_{1};\tilde{A}_{1})
\int_{\tilde{A}_{23}}
Q_{B_{1}B_{2}}^{(\underline{x})}(z_{1};\Gamma_{2})
\tilde{Q}_{B_{1}}(z_{1};d \underline{x}).$$

\noindent The result follows by Lemma \ref{lemmaA2} (Appendix A)
since on the set $B_{1}^{l} \times C^{r-l} \times
(B_{2}^{c})^{n-l}$, $\tilde{Q}_{B_{1}}(z_{1};\cdot)$ is the
product measure between its marginal with respect to the first $l$
components restricted to $B_{1}^{l}$ and its marginal with respect
to the remaining $n-l$ components restricted to $C^{r-l} \times
(B_{2}^{c})^{n-r}$ (by Lemma \ref{Qtilde-disintegration},(a)).
$\Box$

Here is the main result of the paper.

\begin{theorem}
\label{main} If $P:=(P_{A})_{A \in {\cal B}}$ is a ${\cal
Q}$-Markov random  probability measure and $\underline{X}:=(X_{1},
\ldots, X_{n})$ is a sample from $P$, then the family ${\cal
Q}^{(\underline{x})}=(Q_{B_{1}B_{2}}^{(\underline{x})})_{B_{1}
\subseteq B_{2}}$ defined by (\ref{definition-Q-x}) is a posterior
transition system and the conditional distribution of $P$ given
$\underline{X}=\underline{x}$ is ${\cal
Q}^{(\underline{x})}$-Markov.
\end{theorem}

\noindent {\bf Proof}: By Proposition 5 of Balan and Ivanoff
(2002), it is enough to show that $\forall B_{1} \subseteq B_{2}
\subseteq \ldots \subseteq B_{k}$ in ${\cal A}(u)$, $\forall
\tilde{\Gamma} \in {\cal B}([0,1])^{k}$ and for
$\alpha_{n}$-almost all $\underline{x}$
$${\cal P}[(P_{B_{1}}, \ldots,P_{B_{k}}) \in \tilde{\Gamma}|\underline{X}=\underline{x}] =\int_{\tilde{\Gamma}}Q_{B_{k-1}B_{k}}^{(\underline{x})}(z_{k-1};dz_{k}) \ldots
Q_{B_{1}B_{2}}^{(\underline{x})}(z_{1};dz_{2})\mu_{B_{1}}^{(\underline{x})}(dz_{1})$$
or equivalently, for every $\tilde{A} \in {\cal B}^{n}$
\begin{equation}
\label{posterior-chain} {\cal P}(\underline{X} \in \tilde{A},
(P_{B_{j}})_{j} \in \tilde{\Gamma}) = \int_{\tilde{A}}
\int_{\tilde{\Gamma}}Q_{B_{k-1}B_{k}}^{(\underline{x})}(z_{k-1};dz_{k})
\ldots \mu_{B_{1}}^{(\underline{x})}(dz_{1}) \alpha_{n}(d
\underline{x}).
\end{equation}

\noindent Note also that (\ref{posterior-chain}) will imply that
${\cal Q}^{(\underline{x})}$ is a posterior transition system.

For the proof of (\ref{posterior-chain}) we will use an induction
argument on $k \geq 2$. The statement for $k=2$ is exactly
(\ref{disint-nu-B1-B2-n}). Assume that the statement is true for
$k-1$. For each $B_{1} \subseteq B_{2} \subseteq \ldots \subseteq
B_{k}$ in ${\cal A}(u)$ we let $\nu_{B_{1} \ldots B_{k}}$ be the
law of $(X_{1}, \ldots, X_{n}, P_{B_{1}}, \ldots, P_{B_{k}})$
under ${\cal P}$. Note that $\forall A_{1}, \ldots,A_{n} \in {\cal
B}, \forall \Gamma_{1}, \ldots,\Gamma_{k} \in {\cal B}([0,1])$,
$\nu_{B_{1} \ldots B_{k}}(\prod_{i=1}^{n}A_{i} \times
\prod_{j=1}^{k} \Gamma_{j}) = {\cal E}[\prod_{i=1}^{n}P_{A_{i}}
\cdot \prod_{j=1}^{k}I_{\Gamma_{j}}(P_{B_{j}})]$. On the other
hand, $\nu_{B_{1} \ldots B_{k}}(\tilde{A} \times \prod_{j=1}^{k}
\Gamma_{j})$ is also equal to
\begin{equation}
\label{k-sets} \int_{\tilde{A} \times
\prod_{j=1}^{k-1}\Gamma_{j}}Q_{B_{1} \ldots
B_{k}}^{(\underline{x})}(z_{1}, \ldots,z_{k-1};\Gamma_{k})
\nu_{B_{1} \ldots B_{k-1}}(d \underline{x} \times dz_{1} \times
\ldots \times dz_{k-1})=
\end{equation}
$$\int_{\prod_{j=1}^{k}\Gamma_{j}} \tilde{Q}_{B_{1} \ldots B_{k}}(z_{1}, \ldots,z_{k};\tilde{A}) \mu_{B_{1} \ldots B_{k}}(dz_{1} \times \ldots \times dz_{k})$$
where $Q_{B_{1} \ldots B_{k}}^{(\underline{x})}(z_{1},
\ldots,z_{k-1};\Gamma_{k}) := {\cal P}[P_{B_{k}} \in
\Gamma_{k}|\underline{X}=\underline{x}, P_{B_{j}}=z_{j}, j< k]$
and $\tilde{Q}_{B_{1} \ldots B_{k}}(z_{1},
\ldots,z_{k};\prod_{i=1}^{n}A_{i}) := {\cal P}[X_{1} \in A_{1},
 \ldots,X_{n} \in A_{n}|P_{B_{j}}=z_{j}, j \leq k]
 =  {\cal E}[\prod_{i=1}^{n}P_{A_{i}} | P_{B_{j}}=z_{j}, j \leq
 k]$.

\noindent Using the induction hypothesis, the measure $\nu_{B_{1}
\ldots B_{k-1}}$ disintegrates as
$$Q_{B_{k-2}B_{k-1}}^{(\underline{x})}(z_{k-2};dz_{k-1}) \ldots Q_{B_{1}B_{2}}^{(\underline{x})}(z_{1};dz_{2}) \mu_{B_{1}}^{(\underline{x})}(dz_{1}) \alpha_{n}(d \underline{x})$$
Therefore, it is enough to prove that for every $\Gamma_{k} \in
{\cal B}([0,1])$ and for $\nu_{B_{1} \ldots B_{k-1}}$-almost all
$(\underline{x},z_{1}, \ldots,z_{k-1})$
\begin{equation}
\label{Markov-property-n-k} Q_{B_{1} \ldots
B_{k}}^{(\underline{x})}(z_{1},
\ldots,z_{k-1};\Gamma_{k})=Q_{B_{k-1}B_{k}}^{(\underline{x})}(z_{k-1};
\Gamma_{k})
\end{equation}

\noindent On the other hand, the measure $\nu_{B_{1} \ldots
B_{k-1}}$ disintegrates also as
$$\tilde{Q}_{B_{1} \ldots B_{k-1}}(z_{1}, \ldots z_{k-1};d \underline{x}) \mu_{B_{1} \ldots B_{k-1}}(dz_{1} \times \ldots \times dz_{k-1})$$
with respect to its marginal $\mu_{B_{1} \ldots B_{k-1}}$ with
respect to the last $k-1$ components. By the ${\cal Q}$-Markov
property, the measure $\mu_{B_{1} \ldots B_{k}}$ disintegrates as
$$Q_{B_{k-1}B_{k}}(z_{k-1};dz_{k})\mu_{B_{1} \ldots B_{k-1}}(dz_{1} \times \ldots \times dz_{k-1}).$$

\noindent Using (\ref{k-sets}) we can conclude that for
$\mu_{B_{1} \ldots B_{k-1}}$-almost all $(z_{1}, \ldots,z_{k-1})$
\begin{equation}
\label{Bayes-n-k} \int_{\tilde{A}} Q_{B_{1} \ldots
B_{k}}^{(\underline{x})}(z_{1}, \ldots,z_{k-1};\Gamma_{k})
\tilde{Q}_{B_{1} \ldots B_{k-1}}(z_{1}, \ldots,z_{k-1};d
\underline{x})=
\end{equation}
$$\int_{\Gamma_{k}} \tilde{Q}_{B_{1} \ldots B_{k}}(z_{1}, \ldots,z_{k}; \tilde{A}) Q_{B_{k-1}B_{k}}(z_{k-1};dz_{k}).$$

Let $C_{1}=B_{1};C_{j}=B_{j} \verb2\2  B_{j-1},j=2,
\ldots,k;C_{k+1}=B_{k}^{c}$. Note that each $C_{j} \in {\cal
C}(u)$ and $(C_{1}, \ldots,C_{k+1})$ is a partition of ${\cal X}$;
hence each point $x_{i}$ falls into exactly one set of this
partition.

We proceed to the proof of (\ref{Markov-property-n-k}) and we will
suppose that for some $0 \leq l \leq r \leq n$, the points $x_{1},
\ldots,x_{l}$ fall into $B_{k-1}$ (more precisely, each $x_{i}$
falls into some $C_{j_{i}}$ with $1 \leq j_{1} < \ldots < j_{l}
\leq k-1$), the points $x_{l+1}, \ldots,x_{r}$ fall into $C_{k}$
and the points $x_{r+1}, \ldots,x_{n}$ fall into $C_{k+1}$.

The main tool will be (\ref{Bayes-n-k}) where we will consider a
set $\tilde{A}$ of the form $$\tilde{A}:=\prod_{i=1}^{l}(A_{i}
\cap C_{j_{i}}) \times \prod_{i=l+1}^{r}(A_{i} \cap C_{k}) \times
\prod_{i=r+1}^{n}(A_{i} \cap C_{k+1}), \ A_{i} \in {\cal B}.$$

 Let $\tilde{A}_{2}:=\prod_{i=l+1}^{r}(A_{i} \cap C_{k}) \times {\cal X}^{n-r+l}, \tilde{A}_{3}:=\prod_{i=r+1}^{n}(A_{i} \cap C_{k+1}) \times {\cal X}^{r}$ and $\tilde{A}_{23}:=\tilde{A}_{2} \cap \tilde{A}_{3}$. We will prove that
\begin{equation}
\label{Q-tilda-1,k} \tilde{Q}_{B_{1} \ldots B_{k}}(z_{1},
\ldots,z_{k}; \tilde{A})=M \cdot
\tilde{Q}_{B_{k-1}B_{k}}(z_{k-1},z_{k}; \tilde{A}_{23})
\end{equation}
\begin{equation}
\label{Q-tilda-1,k-1} \tilde{Q}_{B_{1} \ldots B_{k-1}}(z_{1},
\ldots,z_{k-1}; \tilde{A})=M \cdot \tilde{Q}_{B_{k-1}}(z_{k-1};
\tilde{A}_{23})
\end{equation}
where $M:=\prod_{i=1}^{l}
\tilde{Q}_{B_{j_{i}-1}B_{j_{i}}}(z_{j_{i}-1},z_{j_{i}}; (A_{i}
\cap C_{j_{i}}) \times {\cal X}^{n-1})$. Then we will have
$$\int_{\tilde{A}} Q_{B_{1} \ldots B_{k}}^{(\underline{x})}(z_{1}, \ldots,z_{k-1};\Gamma_{k}) \tilde{Q}_{B_{1} \ldots B_{k-1}}(z_{1}, \ldots,z_{k-1};d \underline{x})=$$
$$M \cdot \int_{\Gamma_{k}} \tilde{Q}_{B_{k-1}B_{k}}(z_{k-1},z_{k}; \tilde{A}_{23})Q_{B_{k-1}B_{k}}(z_{k-1};dz_{k})=$$
$$M \cdot \int_{\tilde{A}_{23}}Q_{B_{k-1}B_{k}}^{(\underline{x})}(z_{k-1};\Gamma_{k}) \tilde{Q}_{B_{k-1}}(z_{k-1};d \underline{x})=$$
$$\int_{\tilde{A}}Q_{B_{k-1}B_{k}}^{(\underline{x})}(z_{k-1};\Gamma_{k})\tilde{Q}_{B_{1} \ldots B_{k-1}}(z_{1}, \ldots,z_{k-1};d \underline{x})$$
where we used (\ref{Bayes-n-k}) and (\ref{Q-tilda-1,k}) for the
first equality, (\ref{central-Bayes}) for the second equality and
(\ref{Q-tilda-1,k-1}) for the third equality (taking in account
that $Q_{B_{k-1}B_{k}}^{(\underline{x})}(z_{k-1};\Gamma_{k})$ does
not depend on $x_{1}, \ldots,x_{l}$). Relation
(\ref{Markov-property-n-k}) will follow immediately.

It remains to prove (\ref{Q-tilda-1,k}) and (\ref{Q-tilda-1,k-1}).
Using Lemma 3 of Balan and Ivanoff (2002) we have (for $A_{i} \in
{\cal A}$):
$${\cal E}[\prod_{i=r+1}^{n}P_{A_{i} \cap C_{k+1}}|{\cal F}_{B_{k}}] = {\cal E}[\prod_{i=r+1}^{n}P_{A_{i} \cap C_{k+1}}|P_{B_{k}}] = \tilde{Q}_{B_{k}}(P_{B_{k}}; \tilde{A}_{3})$$

\noindent and therefore, by double conditioning with respect to
${\cal F}_{B_{k}}$
$$\tilde{Q}_{B_{1} \ldots B_{k}}((z_{j})_{j \leq k}; \tilde{A})=
{\cal E}[\prod_{i=1}^{l}P_{A_{i} \cap C_{j_{i}}}
\prod_{i=l+1}^{r}P_{A_{i} \cap C_{k}} \prod_{i=r+1}^{n}P_{A_{i}
\cap C_{k+1}}|P_{B_{j}}=z_{j},j \leq k]$$
$$=\tilde{Q}_{B_{k}}(z_{k}; \tilde{A}_{3}) \cdot
{\cal E}[\prod_{i=1}^{l}P_{A_{i} \cap C_{j_{i}}}
\prod_{i=l+1}^{r}P_{A_{i} \cap C_{k}} \ | \ P_{B_{j}}=z_{j},j \leq
k].$$

\noindent For the second term we have
$${\cal E}[\prod_{i=1}^{l}P_{A_{i} \cap C_{j_{i}}} \prod_{i=l+1}^{r}P_{A_{i} \cap C_{k}} \ | \ P_{B_{j}},j \leq k]=$$
$${\cal E}[\prod_{i=1}^{l}P_{A_{i} \cap C_{j_{i}}}  \cdot  {\cal E}[\prod_{i=l+1}^{r}P_{A_{i} \cap C_{k}}|P_{ B_{j_{i}-1} \cup (A_{i} \cap B_{j_{i}}) },i \leq l;P_{B_{j}},j \leq k] \ | \ P_{B_{j}},j \leq k].$$

\noindent Since $P_{A_{i} \cap C_{k}}=P_{B_{k-1} \cup (A_{i} \cap
B_{k})}-P_{B_{k-1}}$, using Lemma \ref{lemmaA1} (Appendix A)
$${\cal E}[\prod_{i=l+1}^{r}P_{A_{i} \cap C_{k}}|P_{B_{j_{i}-1} \cup (A_{i} \cap B_{j_{i}}) },i \leq l;P_{B_{j}}, j \leq k] =\tilde{Q}_{B_{k-1}B_{k}}(P_{B_{k-1}},P_{B_{k}}; \tilde{A}_{2}).$$

\noindent (In order to use Lemma \ref{lemmaA1} we need $A_{l+1}
\subseteq A_{l+2} \subseteq \ldots \subseteq A_{r}$, but this is
not a restriction as we have seen in the proof of Lemma
\ref{Qtilde-disintegration}.)

\noindent Note that by (\ref{Q-tilda-B1-B2''}),
$\tilde{Q}_{B_{k-1}B_{k}}(z_{k-1},z_{k}; \tilde{A}_{2}) \cdot
\tilde{Q}_{B_{k}}(z_{k};
\tilde{A}_{3})=\tilde{Q}_{B_{k-1}B_{k}}(z_{k-1},z_{k};
\tilde{A}_{23})$. Hence the proof of (\ref{Q-tilda-1,k}) will be
complete once we show that
\begin{equation}
\label{tilda-Q-M} {\cal E}[\prod_{i=1}^{l}P_{A_{i} \cap C_{j_{i}}}
\ | \ P_{B_{j}}, j \leq k]=\prod_{i=1}^{l}
\tilde{Q}_{B_{j_{i}-1}B_{j_{i}}}(z_{j_{i}-1},z_{j_{i}}; (A_{i}
\cap C_{j_{i}}) \times {\cal X}^{n-1}).
\end{equation}
But this follows by induction on $l$, using Lemma \ref{lemmaA1}
(Appendix A).

We turn now to the proof of (\ref{Q-tilda-1,k-1}). Using Lemma 3
of Balan and Ivanoff (2002) we have:
$${\cal E}[\prod_{i=l+1}^{r}P_{A_{i} \cap C_{k}} \prod_{i=r+1}^{n}P_{A_{i} \cap C_{k+1}}|{\cal F}_{B_{k-1}}] =\tilde{Q}_{B_{k-1}}(P_{B_{k-1}}; \tilde{A}_{23})$$

\noindent and thereferore, by double conditioning with respect to
${\cal F}_{B_{k-1}}$ we obtain the following expression for
$\tilde{Q}_{B_{1} \ldots B_{k-1}}(z_{1}, \ldots,z_{k-1};
\tilde{A}_{23})$:
$${\cal E}[\prod_{i=1}^{l}P_{A_{i} \cap C_{j_{i}}} \prod_{i=l+1}^{r}P_{A_{i} \cap C_{k}}  \prod_{i=r+1}^{n}P_{A_{i} \cap C_{k+1}}|P_{B_{j}}=z_{j},j \leq k-1]=$$
$$\tilde{Q}_{B_{k-1}}(z_{k-1};  \tilde{A}_{23}) \cdot
{\cal E}[\prod_{i=1}^{l}P_{A_{i} \cap C_{j_{i}}} \ | \
P_{B_{j}}=z_{j},j \leq k-1]$$ and (\ref{Q-tilda-1,k-1}) follows,
using (\ref{tilda-Q-M}). The proof of the theorem is complete.
$\Box$

\vspace{3mm}

The posterior distribution of a Dirichlet process is also
Dirichlet. In the case of an empirical measure which corresponds
to a sample either from a non-random distribution or from a
Dirichlet process, the calculations for the posterior transition
probabilities $Q_{B_{1}B_{2}}^{(\underline{x})}$ are not
straightforward for samples of size greater than $1$; however, in
the case of a sample of size $1$ we have the following result.

\begin{proposition}
If $P:=(P_{A})_{A \in {\cal B}}$ is the empirical measure of a
sample of size $N$ from a non-random distribution $P_{0}$
(respectively from a Dirichlet process with parameter measure
$\alpha$) and $X$ is a sample of size $1$ from $P$, then the
conditional distribution of $P$ given $X=x$ is ${\cal
Q}^{(x)}$-Markov with
$$Q_{B_{1}B_{2}}^{(x)} \left(\frac{m_{1}}{N}; \left\{ \frac{m_{2}}{N} \right\}\right)=
Q_{B_{1}B_{2}}^{(1)} \left(\frac{m_{1}- \delta_{x}(B_{1})}{N-1};
\left\{ \frac{m_{2}- \delta_{x}(B_{2})}{N-1} \right\}\right)$$
where ${\cal Q}^{(1)}$ is the transition system of the empirical
measure of a sample of size $N-1$ from $P_{0}$ (respectively from
a Dirichlet process with parameter measure $\alpha$).
\end{proposition}

\noindent {\bf Proof}: Let $P$ be the empirical measure of a
sample from a non-random distribution $P_{0}$. Note that
$\alpha_{1}(A)={\cal E}[P_{A}]=P_{0}(A), \forall A \in {\cal B}$.
We have
$$Q_{B_{1}B_{2}}^{(x)} \left(\frac{m_{1}}{N}; \left\{ \frac{m_{2}}{N} \right\}\right)=
Q_{B_{1}B_{2}} \left(\frac{m_{1}}{N}; \left\{ \frac{m_{2}}{N}
\right\}\right)= Q_{B_{1}B_{2}}^{(1)} \left(\frac{m_{1}-1}{N-1};
\left\{ \frac{m_{2}-1}{N-1} \right\}\right)$$ for
$\alpha_{1}$-almost all $x \in B_{1}$. The fact that
$$Q_{B_{1}B_{2}}^{(x)} \left(\frac{m_{1}}{N}; \left\{ \frac{m_{2}}{N} \right\}\right)=
Q_{B_{1}B_{2}}^{(1)} \left(\frac{m_{1}}{N-1}; \left\{
\frac{m_{2}-1}{N-1} \right\}\right)$$ for $\alpha_{1}$-almost all
$x \in C$ follows from (\ref{central-Bayes}),
since for every $A \in {\cal B}$
$$\tilde{Q}_{B_{1}}\left( \frac{m_{1}}{N};A \cap C \right)={\cal E}[P_{A \cap C}|P_{B_{1}}=\frac{m_{1}}{N}]=\frac{N-m_{1}}{N} \cdot \frac{P_{0}(A \cap C)}{P_{0}(B_{1}^{c})}$$
$$\tilde{Q}_{B_{1}B_{2}}\left( \frac{m_{1}}{N},\frac{m_{2}}{N};A \cap C \right)={\cal E}[P_{A \cap C}|P_{B_{1}}=\frac{m_{1}}{N},P_{B_{2}}=\frac{m_{2}}{N}]=\frac{m_{2}-m_{1}}{N} \cdot \frac{P_{0}(A \cap C)}{P_{0}(C)}.$$

\noindent Similarly one can show that
$$\tilde{Q}_{B_{1}}\left( \frac{m_{1}}{N};A \cap B_{2}^{c} \right)=\frac{N-m_{1}}{N} \cdot \frac{P_{0}(A \cap B_{2}^{c})}{P_{0}(B_{1}^{c})}$$
$$\tilde{Q}_{B_{1}B_{2}}\left( \frac{m_{1}}{N},\frac{m_{2}}{N};A \cap B_{2}^{c} \right)=\frac{N-m_{2}}{N} \cdot \frac{P_{0}(A \cap B_{2}^{c})}{P_{0}(B_{2}^{c})}$$
and hence
$$Q_{B_{1}B_{2}}^{(x)} \left(\frac{m_{1}}{N}; \left\{ \frac{m_{2}}{N} \right\}\right)=
Q_{B_{1}B_{2}}^{(1)} \left(\frac{m_{1}}{N-1}; \left\{
\frac{m_{2}}{N-1} \right\}\right)$$ for $\alpha_{1}$-almost all
$x$ in $B_{2}^{c}$.

If $P$ is the empirical measure of a sample from a Dirichlet
process with parameter measure $\alpha$, then
$\alpha_{1}(A)=\alpha(A)/\alpha({\cal X})$ and a similar argument
can be used. $\Box$

\section{Neutral to the right random probability measures}

Let $P:=(P_{A})_{A \in {\cal B}}$ be a random probability measure
on ${\cal X}$. For every sets $B_{1},B_{2} \in {\cal A}(u)$ with
$B_{1} \subseteq B_{2}$, we define $V_{B_{1}B_{2}}$ to be equal to
$(P_{B_{2}}-P_{B_{1}})/(1-P_{B_{1}})$ on the set $\{P_{B_{1}}<1\}$
and $1$ elsewhere; let $F_{B_{1}B_{2}}$ be the distribution of
$V_{B_{1}B_{2}}$. The next definition generalizes the definition
of Doksum (1974).

\begin{definition}
{\rm A random probability measure $P:=(P_{A})_{A \in {\cal B}}$ is
called {\em neutral to the right} if for every sets $B_{1}
\subseteq \ldots \subseteq B_{k}$ in ${\cal A}(u)$, $P_{B_{1}},
V_{B_{1}B_{2}}, \ldots,V_{B_{k-1}B_{k}}$ are independent.}
\end{definition}

{\em Comments}: 1. A random probability measure $P:=(P_{A})_{A \in
{\cal B}}$ is neutral to the right if and only if $\forall
B_{1},B_{2} \in {\cal A}(u),B_{1} \subseteq B_{2}$,
$V_{B_{1}B_{2}}$ is independent of ${\cal F}_{B_{1}}$.

2. The Dirichlet process with parameter measure $\alpha$ is
neutral to the right with $F_{B_{1}B_{2}}$ equal to the Beta
distribution with parameters $\alpha(B_{2} \verb2\2 B_{1}),
\alpha(B_{2}^{c})$.

3. If we denote $C_{1}=B_{1}; C_{i}=B_{i} \verb2\2 B_{i-1};i=2,
\ldots,k$, then $(P_{C_{1}}, \ldots,P_{C_{k}})$ has a `completely
neutral' distribution (see Definition \ref{completely-neutral});
this distribution was formally introduced by Connor and Mosimann
(1969), although the concept itself goes back to Halmos (1944).
Note that the Dirichlet process is the only non-trivial process
which has completely neutral distributions over any disjoint sets
$\{A_{1}, \ldots,A_{k}\}$ in ${\cal B}$ (according to Ferguson
1974, p. 622).

4. In general, the process $Y_{A}:=- \ln(1-P_{A}), A \in {\cal B}$
is not additive and hence it does not have independent increments,
even if $Y_{B_{1}}, Y_{B_{2}}-Y_{B_{1}}, \ldots,
Y_{B_{k}}-Y_{B_{k-1}}$ are independent for any sets $B_{1}
\subseteq B_{2} \subseteq \ldots \subseteq B_{k}$ in ${\cal A}(u)$
(the increment $Y_{B_{2} \verb2\2 B_{1}}$ is not equal to
$Y_{B_{2}}-Y_{B_{1}}$); therefore, the theory of processes with
independent increments cannot be used in higher dimensions.

\begin{proposition}
\label{neutral-QMarkov} A neutral to the right random probability
measure is ${\cal Q}$-Markov with
\begin{equation}
\label{trans-neutral} Q_{B_{1}B_{2}}(z_{1}; \Gamma_{2}):=
 \left\{
     \begin{array}{ll}
       F_{B_{1}B_{2}} \left( \frac{\Gamma_{2}-z_{1}}{1-z_{1}} \right) & \mbox{if $z_{1}<1$} \\
       \delta_{1}(\Gamma_{2})  & \mbox{if $z_{1}=1$}
     \end{array}
 \right.
\end{equation}
\end{proposition}

\noindent {\bf Proof}: For any sets $B_{1} \subseteq \ldots
\subseteq B_{k}$ in ${\cal A}(u)$, $P_{B_{1}}, \ldots, P_{B_{k}}$
is a Markov chain:
\begin{eqnarray*}
\lefteqn{ {\cal P}[P_{B_{j}} \in \Gamma_{j}|P_{B_{1}}=z_{1}, \ldots, P_{B_{j-1}}=z_{j-1}] = } \\
 & & {\cal P} [ V_{B_{j-1}B_{j}} \in \frac{\Gamma_{j}-z_{j-1}}{1-z_{j-1}}|P_{B_{1}}=z_{1}, V_{B_{1}B_{2}}=v_{2}, \ldots, V_{B_{j-2}B_{j-1}}=v_{j-1} ] = \\
 & & {\cal P} [ V_{B_{j-1}B_{j}} \in \frac{\Gamma_{j}-z_{j-1}}{1-z_{j-1}} ] =
{\cal P} [ V_{B_{j-1}B_{j}} \in \frac{\Gamma_{j}-z_{j-1}}{1-z_{j-1}}|P_{B_{j-1}}=z_{j-1}] = \\
 & & {\cal P}[P_{B_{j}} \in \Gamma_{j}|P_{B_{j-1}}=z_{j-1}]
\end{eqnarray*}
where $v_{i}:=(z_{i}-z_{i-1})/(1-z_{i-1}),i=2, \ldots,j-1$ and
assuming $z_{i}<1, \forall i$. $\Box$

\vspace{2mm}

For any sets $B_{1} \subseteq B_{2} \subseteq B_{3}$ in ${\cal
A}(u)$,
$V_{B_{1}B_{3}}=V_{B_{1}B_{2}}+V_{B_{2}B_{3}}-V_{B_{1}B_{2}} \cdot
V_{B_{2}B_{3}}$. This leads us to the following definition.

\begin{definition}
\label{definition-F} For each $B_{1}, B_{2} \in {\cal A}(u)$ with
$B_{1} \subseteq B_{2}$, let $F_{B_{1}B_{2}}$ be a probability
measure on $[0,1]$. The family $(F_{B_{1}B_{2}})_{B_{1} \subseteq
B_{2}}$ is called a {\bf neutral to the right system} if $\forall
B_{1} \subseteq B_{2} \subseteq B_{3}$ in ${\cal A}(u)$
$$F_{B_{1}B_{3}}(\Gamma)=\int_{[0,1]^{2}} I_{\Gamma}(y+z-yz) F_{B_{2}B_{3}}(dz) F_{B_{1}B_{2}}(dy).$$
\end{definition}

{\em Comments}: 1. If we let $U_{B_{1}B_{2}}:=-
\ln(1-V_{B_{1}B_{2}})$ and $G_{B_{1}B_{2}}$ be the distribution of
$U_{B_{1}B_{2}}$, then for every $B_{1} \subseteq B_{2} \subseteq
B_{3}$ in ${\cal A}(u)$,
$U_{B_{1}B_{3}}=U_{B_{1}B_{2}}+U_{B_{2}B_{3}}$ and
$G_{B_{1}B_{3}}=G_{B_{1}B_{2}} \ast G_{B_{2}B_{3}}$.

2. Let $Q_{B_{1}B_{2}}(z_{1};
\Gamma_{2}):=F_{B_{1}B_{2}}((\Gamma_{2}-z_{1})/(1-z_{1}))$ for
$z_{1}<1$ and $Q_{B_{1}B_{2}}(1; \cdot)=\delta_{1}$; then
$(F_{B_{1}B_{2}})_{B_{1} \subseteq B_{2}}$ is a neutral to the
right system if and only if $(Q_{B_{1}B_{2}})_{B_{1} \subseteq
B_{2}}$ is a transition system.

\vspace{2mm}

The following result is the converse of Proposition
\ref{neutral-QMarkov}.

\begin{proposition}
\label{QMarkov-neutral} If $P:=(P_{A})_{A \in {\cal B}}$ is a
${\cal Q}$-Markov random probability measure with a transition
system ${\cal Q}$ given by (\ref{trans-neutral}) for a neutral to
the right system $(F_{B_{1}B_{2}})_{B_{1} \subseteq B_{2}}$, then
$P$ is neutral to the right.
\end{proposition}

\noindent {\bf Proof}: We want to prove that for every $B_{1},
B_{2} \in {\cal A}(u)$ with $B_{1} \subseteq B_{2}$ and for every
$A_{1}, \ldots, A_{k} \in {\cal A}, A_{i} \subseteq
B,A_{k}=B_{1}$, $V_{B_{1}B_{2}}$ is independent of $(P_{A_{1}},
\ldots,P_{A_{k}})$. Using the ${\cal Q}$-Markov property we have:
${\cal P}[V_{B_{1}B_{2}} \in \Gamma|P_{A_{i}}=z_{i}; i=1,\ldots,k]
= {\cal P}[P_{B_{2}} \in z_{k}+(1-z_{k})\Gamma|P_{B_{1}}=z_{k}]=
Q_{B_{1}B_{2}}(z_{k}; z_{k}+(1-z_{k})\Gamma)=
F_{B_{1}B_{2}}(\Gamma)={\cal P}(V_{B_{1}B_{2}} \in \Gamma)$. Since
this holds for any Borel set $\Gamma$ in $[0,1]$, the proof is
complete. $\Box$

In what follows we will prove that the posterior distribution of a
neutral to the right random probability measure is also neutral to
the right, by showing that the posterior transition probabilities
$Q_{B_{1}B_{2}}^{(\underline{x})}$ are of the form
(\ref{trans-neutral}) for a ``posterior'' neutral to the right
system $(F_{B_{1}B_{2}}^{(\underline{x})})_{B_{1} \subseteq
B_{2}}$. This extends Doksum's (1974) result to an arbitrary space
${\cal X}$, which can be endowed with an indexing collection
${\cal A}$.

Let $P:=(P_{A})_{A \in {\cal B}}$ be a neutral to the right
process and $\underline{X}:=(X_{1}, \ldots, X_{n})$ a sample from
$P$. In order to define the probability measures
$F_{B_{1}B_{2}}^{(\underline{x})}$ we will use the same Bayesian
technique as in Section 3.

For each sets $B_{1},B_{2} \in {\cal A}(u);B_{1} \subseteq B_{2}$,
let $\phi_{B_{1}B_{2}}$ be the law of $X_{1}, \ldots,X_{n},
\linebreak V_{B_{1}B_{2}}$ under ${\cal P}$. Note that
$\phi_{B_{1}B_{2}}(\prod_{i=1}^{n}A_{i} \times \Gamma)={\cal
E}[\prod_{i=1}^{n}P_{A_{i}} \cdot I_{\Gamma}(V_{B_{1}B_{2}})]$. On
the other hand, we have
\begin{equation}
\label{disint-phi} \phi_{B_{1}B_{2}}(\tilde{A} \times \Gamma)=
\int_{\tilde{A}} F_{B_{1}B_{2}}^{(\underline{x})}(\Gamma)
\alpha_{n}(d \underline{x})
 =  \int_{\Gamma} \tilde{T}_{B_{1}B_{2}}(z; \tilde{A}) F_{B_{1}B_{2}}(dz)
\end{equation}
where
\begin{equation}
\label{definition-F-x}
F_{B_{1}B_{2}}^{(\underline{x})}(\Gamma):={\cal P}[V_{B_{1}B_{2}}
\in \Gamma|\underline{X}= \underline{x}]
\end{equation}
and $\tilde{T}_{B_{1}B_{2}}(z; \tilde{A}):={\cal P}[\underline{X}
\in \tilde{A}|V_{B_{1}B_{2}}=z]$.

In the proof of Theorem \ref{main-neutral} we will see that
$(F_{B_{1}B_{2}}^{(\underline{x})})_{B_{1} \subseteq B_{2}}$ may
not be a genuine neutral to the right system as introduced by
Definition \ref{definition-F}. Therefore we need to introduce the
following terminology.

\begin{definition}
The family $(F_{B_{1}B_{2}}^{(\underline{x})})_{B_{1} \subseteq
B_{2}}$ defined by (\ref{definition-F-x}) is called a {\bf
posterior neutral to the right system} (corresponding to $P$ and
$\underline{X}$) if $\forall B_{1} \subseteq B_{2} \subseteq
B_{3}$ in ${\cal A}(u)$, $\forall \Gamma \in {\cal B}([0,1])$ and
for $\alpha_{n}$-almost all $\underline{x}$
$$F_{B_{1}B_{3}}^{(\underline{x})}(\Gamma)= \int_{[0,1]^{2}}
I_{\Gamma}(y+z-yz)
F_{B_{2}B_{3}}^{(\underline{x})}(dz)F_{B_{1}B_{2}}^{(\underline{x})}(dy).$$

\noindent The conditional distribution of $P$ given
$\underline{X}=\underline{x}$ is called {\bf neutral to the right}
if $\forall B_{1} \subseteq B_{2}$ in ${\cal A}(u)$,
$V_{B_{1}B_{2}}$ is conditionally independent of ${\cal
F}_{B_{1}}$ given $\underline{X}$.
\end{definition}

 Let $C:=B_{2} \verb2\2 B_{1}$. For fixed $0 \leq l \leq r
\leq n$ we will consider sets of the form
$\tilde{A}_{23}:=\prod_{i=l+1}^{r}(A_{i} \cap C) \times
\prod_{i=r+1}^{n}(A_{i} \cap B_{2}^{c})$, where $A_{i} \in {\cal
B}$.

\begin{lemma}
\label{tildeQ-neutral-23} (a) For $\mu_{B_{1}}$-almost all
$z_{1}$,
$$\tilde{Q}_{B_{1}}(z_{1}; \tilde{A}_{23}) = \frac{(1-z_{1})^{n-l}}{\alpha_{n}((B_{1}^{c})^{n-l} \times {\cal X}^{l})} \cdot \alpha_{n}(\tilde{A}_{23}).$$

(b) For $\mu_{B_{1}B_{2}}$-almost all $(z_{1},z_{2})$,
$$\tilde{Q}_{B_{1}B_{2}}(z_{1},z_{2}; \tilde{A}_{23}) = \frac{(1-z_{1})^{n-l}}{\alpha_{n}((B_{1}^{c})^{n-l} \times {\cal X}^{l})} \cdot \tilde{T}_{B_{1}B_{2}}
\left( \frac{z_{2}-z_{1}}{1-z_{1}};\tilde{A}_{23}\right).$$
\end{lemma}

\noindent {\bf Proof}: Without loss of generality we will assume
that $A_{i} \in {\cal A}, \forall i$. We have
\begin{equation}
\label{neutral-Q-tilde} \prod_{i=l+1}^{r}P_{A_{i} \cap C} \cdot
\prod_{i=r+1}^{n}P_{A_{i} \cap B_{2}^{c}} = (1-P_{B_{1}})^{n-l}
\cdot \prod_{i=l+1}^{r} \frac{P_{A_{i} \cap C}}{1-P_{B_{1}}} \cdot
\prod_{i=r+1}^{n} \frac{P_{A_{i} \cap B_{2}^{c}}}{1-P_{B_{1}}}.
\end{equation}

\noindent Note that $P_{A_{i} \cap
C}/(1-P_{B_{1}})=V_{B_{1},(A_{i} \cap B_{2}) \cup B_{1}}, P_{A_{i}
\cap B_{2}^{c}}/(1-P_{B_{1}})=V_{B_{1}, A_{i} \cup
B_{2}}-V_{B_{1}B_{2}}$ and $P_{B_{1}}$ is independent of
$V_{B_{1}, (A_{i} \cap B_{2}) \cup B_{1}}, i=l+1, \ldots, r,
V_{B_{1}B_{2}}$ and $V_{B_{1},A_{i} \cup B_{2}}, i=r+1, \ldots,
n$.

\noindent (a) Take ${\cal E}[ \ \cdot \ ]$, respectively ${\cal
E}[\ \cdot \ |P_{B_{1}}=z_{1}]$ in (\ref{neutral-Q-tilde}); we get
\begin{equation}
\label{expect-tildeA23} {\cal E}[\prod_{i=l+1}^{r} \frac{P_{A_{i}
\cap C}}{1-P_{B_{1}}} \cdot \prod_{i=r+1}^{n} \frac{P_{A_{i} \cap
B_{2}^{c}}}{1-P_{B_{1}}}]=
\frac{\alpha_{n}(\tilde{A}_{23})}{\alpha_{n}((B_{1}^{c})^{n-l}
\times {\cal X}^{l})}
\end{equation}
$$\tilde{Q}_{B_{1}}(z_{1}; \tilde{A}_{23})=(1-z_{1})^{n-l} \cdot {\cal E} [\prod_{i=l+1}^{r}P_{A_{i} \cap C}  \cdot \prod_{i=r+1}^{n}P_{A_{i} \cap B_{2}^{c}}]=\frac{(1-z_{1})^{n-l} \alpha_{n}(\tilde{A}_{23})}{\alpha_{n}((B_{1}^{c})^{n-l} \times {\cal X}^{l})}.$$

\noindent (b) 
Take ${\cal E}[\ \cdot \ |V_{B_{1}B_{2}}=z]$, respectively ${\cal
E}[\ \cdot \ |P_{B_{1}}=z_{1},P_{B_{2}}=z_{2}]$ in
(\ref{neutral-Q-tilde}); we get
\begin{equation}
\label{cond-expect-tildeA23} {\cal
E}[\prod_{i=l+1}^{r}\frac{P_{A_{i} \cap C}}{1-P_{B_{1}}}
\prod_{i=r+1}^{n} \frac{P_{A_{i} \cap
B_{2}^{c}}}{1-P_{B_{1}}}|V_{B_{1}B_{2}}=z]=\frac{\tilde{T}_{B_{1}B_{2}}(z;\tilde{A}_{23})}{\alpha_{n}((B_{1}^{c})^{n-l}
\times {\cal X}^{l})}
\end{equation}
$$\tilde{Q}_{B_{1}B_{2}}(z_{1},z_{2}; \tilde{A}_{23})=(1-z_{1})^{n-l}{\cal E} [\prod_{i=l+1}^{r}\frac{P_{A_{i} \cap C}}{1-P_{B_{1}}} \prod_{i=r+1}^{n}\frac{P_{A_{i} \cap B_{2}^{c}}}{1-P_{B_{1}}}|V_{B_{1}B_{2}}=\frac{z_{2}-z_{1}}{1-z_{1}}]$$
$$=(1-z_{1})^{n-l} \cdot \frac{1}{\alpha_{n}((B_{1}^{c})^{n-l} \times {\cal X}^{l})} \cdot \tilde{T}_{B_{1}B_{2}} \left(\frac{z_{2}-z_{1}}{1-z_{1}};\tilde{A}_{23} \right)$$
which concludes the proof. $\Box$

\begin{lemma}
For every $B_{1},B_{2} \in {\cal A}(u)$ with $B_{1} \subseteq
B_{2}$, for every $\Gamma \in {\cal B}([0,1])$ and for
$\alpha_{n}$-almost all $\underline{x}$,
$F_{B_{1}B_{2}}^{(\underline{x})}(\Gamma)$ does not depend on
those $x_{i}$'s that fall in $B_{1}$; in particular, for
$\alpha_{n}$-almost all $\underline{x}$ in $B_{1}^{n}$,
$F_{B_{1}B_{2}}^{(\underline{x})}(\Gamma)=F_{B_{1}B_{2}}(\Gamma)$.
\end{lemma}

\noindent {\bf Proof}: For arbitrary $A_{1}, \ldots, A_{n} \in
{\cal A}$ we write
\begin{eqnarray*}
\lefteqn{\prod_{i=1}^{l}P_{A_{i} \cap B_{1}} \prod_{i=l+1}^{r}P_{A_{i} \cap C}  \prod_{i=r+1}^{n}P_{A_{i} \cap B_{2}^{c}} = } \\
 & & (1-P_{B_{1}})^{n-l} \prod_{i=1}^{l}P_{A_{i} \cap B_{1}}  \prod_{i=l+1}^{r} \frac{P_{A_{i} \cap C}}{1-P_{B_{1}}} \prod_{i=r+1}^{n} \frac{P_{A_{i} \cap B_{2}^{c}}}{1-P_{B_{1}}}.
\end{eqnarray*}

\noindent Taking ${\cal E}[\ \cdot \ ], {\cal E}[\ \cdot \
|V_{B_{1}B_{2}}=z]$ and using (\ref{expect-tildeA23}),
respectively (\ref{cond-expect-tildeA23}) we get
$$\alpha_{n}(\tilde{A})=\frac{\alpha_{n}(\prod_{i=1}^{l}(A_{i} \cap B_{1}) \times (B_{1}^{c})^{n-l})}{\alpha_{n}((B_{1}^{c})^{n-l} \times {\cal X}^{l})} \cdot \alpha_{n}(\tilde{A}_{23})$$
$$\tilde{T}_{B_{1}B_{2}}(z; \tilde{A})=\frac{\alpha_{n}(\prod_{i=1}^{l}(A_{i} \cap B_{1}) \times (B_{1}^{c})^{n-l})}{\alpha_{n}((B_{1}^{c})^{n-l} \times {\cal X}^{l})} \cdot \tilde{T}_{B_{1}B_{2}}(z; \tilde{A}_{23}).$$

\noindent Using (\ref{disint-phi}) we get
$$\int_{\tilde{A}} F_{B_{1}B_{2}}^{(\underline{x})}(\Gamma) \alpha_{n}(d \underline{x})
= \int_{\Gamma} \tilde{T}_{B_{1}B_{2}}(z; \tilde{A})
F_{B_{1}B_{2}}(dz)=$$
$$\frac{\alpha_{n}(\prod_{i=1}^{l}(A_{i} \cap B_{1}) \times (B_{1}^{c})^{n-l})}{\alpha_{n}((B_{1}^{c})^{n-l} \times {\cal X}^{l})} \cdot  \int_{\Gamma} \tilde{T}_{B_{1}B_{2}}(z; \tilde{A}_{23}) F_{B_{1}B_{2}}(dz)=$$
$$\frac{\alpha_{n}(\prod_{i=1}^{l}(A_{i} \cap B_{1}) \times (B_{1}^{c})^{n-l})}{\alpha_{n}((B_{1}^{c})^{n-l} \times {\cal X}^{l})} \cdot \int_{\tilde{A}_{23}} F_{B_{1}B_{2}}^{(\underline{x})}(\Gamma) \alpha_{n}(d \underline{x}).$$

\noindent The result follows by Lemma \ref{lemmaA2} (Appendix A).
$\Box$

\vspace{3mm}

Here is the main result of this section.

\begin{theorem}
\label{main-neutral} If $P:=(P_{A})_{A \in {\cal B}}$ is a neutral
to the right random probability measure and
$\underline{X}:=(X_{1}, \ldots,X_{n})$ is a sample from $P$, then
the conditional distribution of $P$ given
$\underline{X}=\underline{x}$ is also neutral to the right.
\end{theorem}

\noindent {\bf Proof}: Since $P$ is ${\cal Q}$-Markov, by Theorem
\ref{main} the conditional distribution of $P$ given
$\underline{X}=\underline{x}$ is ${\cal
Q}^{(\underline{x})}$-Markov. Using Lemma \ref{tildeQ-neutral-23},
the key equation (\ref{central-Bayes}) becomes
$$\int_{\tilde{A}_{23}}{\cal Q}_{B_{1}B_{2}}^{(\underline{x})}(z_{1}; \Gamma_{2}) \alpha_{n}(d \underline{x})= \int_{\Gamma_{2}} \tilde{T}_{B_{1}B_{2}} \left( \frac{z_{2}-z_{1}}{1-z_{1}}; \tilde{A}_{23}\right) Q_{B_{1}B_{2}}(z_{1};dz_{2}).$$

\noindent Using Proposition \ref{neutral-QMarkov} and relation
(\ref{disint-phi}), the right-hand side becomes (for $z_{1}<1$)
$$\int_{\frac{\Gamma_{2}-z_{1}}{1-z_{1}}} \tilde{T}_{B_{1}B_{2}}(z; \tilde{A}_{23})F_{B_{1}B_{2}}(dz)=\int_{\tilde{A}_{23}}F_{B_{1}B_{2}}^{(\underline{x})} \left(\frac{\Gamma_{2}-z_{1}}{1-z_{1}} \right) \alpha_{n}(d \underline{x}).$$

\noindent This proves that $\forall z_{1} \in [0,1), \forall
\Gamma_{2} \in {\cal B}([0,1])$ and for $\alpha_{n}$-almost all
$\underline{x}$
$${\cal Q}_{B_{1}B_{2}}^{(\underline{x})}(z_{1}; \Gamma_{2})=
F_{B_{1}B_{2}}^{(\underline{x})}
\left(\frac{\Gamma_{2}-z_{1}}{1-z_{1}} \right).$$

\noindent Since ${\cal Q}^{(\underline{x})}$ is a posterior
transition system, it follows that
$(F_{B_{1}B_{2}}^{(\underline{x})})_{B_{1} \subseteq B_{2}}$ is a
posterior neutral to the right system. By Proposition
\ref{QMarkov-neutral}, the distribution of $P$ given
$\underline{X}$ is neutral to the right. $\Box$

\vspace{2mm}

The next result gives some simple formulas for calculating the
posterior distribution of $P_{B_{1}}$ when all the observations
fall outside $B_{1}$, and the posterior distribution of
$V_{B_{1}B_{2}}$ when all the observations fall outside $B_{2}
\verb2\2 B_{1}$.

\begin{proposition}
(a) For $\alpha_{n}$-almost all $\underline{x}$ with $x_{i} \in
B_{1}^{c} \ \forall i$
$$\mu_{B_{1}}^{(\underline{x})}(\Gamma)={\cal P}[P_{B_{1}} \in \Gamma |\underline{X}=\underline{x}]= \frac{{\cal E}[I_{\Gamma}(P_{B_{1}})(1-P_{B_{1}})^{n}]}{{\cal E}[(1-P_{B_{1}})^{n}]}.$$

(b) For $\alpha_{n}$-almost all $\underline{x}$ with $x_{i} \in
(B_{2} \verb2\2 B_{1})^{c} \ \forall i$
$$F_{B_{1}B_{2}}^{(\underline{x})}(\Gamma)={\cal P}[V_{B_{1}B_{2}} \in \Gamma|\underline{X}=\underline{x}]=\frac{{\cal E}[I_{\Gamma}(V_{B_{1}B_{2}})(1-P_{B_{2}})^{m}]}{{\cal E}[(1-P_{B_{2}})^{m}]}$$
where $m$ denotes the number of $x_{i}$'s that fall outside
$B_{2}$.
\end{proposition}

\noindent {\bf Proof}: Note that (a) is a particular case of (b)
since $\mu_{B_{1}}^{(\underline{x})}=F_{\emptyset
B_{1}}^{(\underline{x})}$. We proceed to the proof of (b). For
fixed $0 \leq l \leq n$, let $\tilde{A}:=\prod_{i=1}^{l}(A_{i}
\cap B_{1}) \times \prod_{i=l+1}^{n}(A_{i} \cap B_{2}^{c})$, where
$A_{i} \in {\cal B}$. We claim that
\begin{equation}
\label{tildeT-outsideC}
\tilde{T}_{B_{1}B_{2}}(z;\tilde{A})=(1-z)^{n-l} \cdot
\frac{\alpha_{n}((B_{1}^{c})^{n-l} \times {\cal
X}^{l})}{\alpha_{n}((B_{2}^{c})^{n-l} \times {\cal X}^{l})} \cdot
\alpha_{n}(\tilde{A})
\end{equation}

\noindent Using (\ref{disint-phi}), it follows that
$$\int_{\tilde{A}} F_{B_{1}B_{2}}^{(\underline{x})}(\Gamma) \alpha_{n}(d \underline{x})
 =  \frac{\alpha_{n}((B_{1}^{c})^{n-l} \times {\cal X}^{l})}{\alpha_{n}((B_{2}^{c})^{n-l} \times {\cal X}^{l})} \cdot \alpha_{n}(\tilde{A}) \cdot
\int_{\Gamma} (1-z)^{n-l} F_{B_{1}B_{2}}(dz)$$ and hence for
$\alpha_{n}$-almost all $\underline{x}$ with $x_{i} \in (B_{2}
\verb2\2 B_{1})^{c}, \forall i$
$$F_{B_{1}B_{2}}^{(\underline{x})}(\Gamma)=\frac{\alpha_{n}((B_{1}^{c})^{n-l} \times {\cal X}^{l})}{\alpha_{n}((B_{2}^{c})^{n-l} \times {\cal X}^{l})} \cdot  \int_{\Gamma} (1-z)^{n-l} F_{B_{1}B_{2}}(dz)=$$
$$\frac{{\cal E}[(1-P_{B_{1}})^{n-l}]}{{\cal E}[(1-P_{B_{1}})^{n}]} \cdot {\cal E}[I_{\Gamma}(V_{B_{1}B_{2}}) (1-V_{B_{1}B_{2}})^{n-l}]=
\frac{{\cal E}[I_{\Gamma}(V_{B_{1}B_{2}})
(1-P_{B_{2}})^{n-l}]}{{\cal E}[(1-P_{B_{1}})^{n}]}$$ since
$P_{B_{1}}$ is independent of $V_{B_{1}B_{2}}$.

We turn now to the proof of (\ref{tildeT-outsideC}). Without loss
of generality we will assume that $A_{i} \in {\cal A}, \forall i$.
Let $\tilde{A}_{2}= \prod_{i=l+1}^{n}(A_{i} \cap B_{2}^{c}) \times
{\cal X}^{l}$. We have
$$\prod_{i=1}^{l}P_{A_{i} \cap B_{1}} \prod_{i=l+1}^{n}P_{A_{i} \cap B_{2}^{c}}=
(1-P_{B_{1}})^{n-l} \prod_{i=1}^{l} P_{A_{i} \cap B_{1}}
(1-V_{B_{1}B_{2}})^{n-l} \prod_{i=l+1}^{n} \frac{P_{A_{i} \cap
B_{2}^{c}}}{1-P_{B_{2}}}.$$

\noindent Note that $(1-P_{B_{1}})^{n-l} \prod_{i=1}^{l} P_{A_{i}
\cap B_{1}}$ is ${\cal F}_{B_{1}}$-measurable and ${\cal
F}_{B_{1}}$ is independent of $V_{B_{1}B_{2}},V_{B_{2},A_{i} \cup
B_{2}},i=l+1, \ldots,n$. By taking ${\cal E}[\ \cdot \
|V_{B_{1}B_{2}}=z]$ we get
$$\tilde{T}_{B_{1}B_{2}}(z;\tilde{A})=(1-z)^{n-l} \cdot {\cal E}[\prod_{i=1}^{l} P_{A_{i} \cap B_{1}} \cdot (1-P_{B_{1}})^{n-l}] \cdot {\cal E}[\prod_{i=l+1}^{n}\frac{P_{A_{i} \cap B_{2}^{c}}}{1-P_{B_{2}}}]=$$
$$(1-z)^{n-l} \cdot \alpha_{n}(\prod_{i=1}^{l}(A_{i} \cap B_{1}) \times (B_{1}^{c})^{n-l}) \cdot \frac{\alpha_{n}(\tilde{A}_{2})}{\alpha_{n}((B_{2}^{c})^{n-l} \times {\cal X}^{l})}.$$

\noindent Finally, by taking expectation in
$$\prod_{i=1}^{l}P_{A_{i} \cap B_{1}} \cdot \prod_{i=l+1}^{n}P_{A_{i} \cap B_{2}^{c}}= \prod_{i=1}^{l} P_{A_{i} \cap B_{1}} \cdot (1-P_{B_{1}})^{n-l} \cdot \prod_{i=l+1}^{n} \frac{P_{A_{i} \cap B_{2}^{c}}}{1-P_{B_{1}}}$$
we get $\alpha_{n}(\tilde{A}) =\alpha_{n}(\prod_{i=1}^{l}(A_{i}
\cap B_{1}) \times (B_{1}^{c})^{n-l}) \cdot
\alpha_{n}(\tilde{A}_{2})/\alpha_{n}((B_{1}^{c})^{n-l} \times
{\cal X}^{l})$. $\Box$

\vspace{3mm}

\footnotesize{{\em Acknowledgement.} I would like to thank
Jyotirmoy Dey and Giovanni Petris for drawing my attention to the
original works of Ferguson and Doksum. I would also like to thank
Professor Jean Vaillancourt for the useful discussions which have
helped clarify the ideas. Finally, I am very grateful to an
anonymous referee who read the paper very carefully and made
numerous suggestions for improvement.

\normalsize{

\appendix

\section{Some elementary results}

\begin{lemma}
\label{lemmaA1} If $(X_{t})_{t \in {\bf R}}$ is a Markov process,
then for every $s_{1} < \ldots <s_{n}<s<u_{1}< \ldots <u_{p}< t<
t_{1}< \ldots< t_{m}$ and for every bounded measurable function
$h$
$${\cal E}[h(X_{u_{1}}, \ldots,X_{u_{p}})|X_{s_{1}}, \ldots,X_{s_{n}},X_{s},X_{t},X_{t_{1}}, \ldots,X_{t_{m}}]={\cal E}[h(X_{u_{1}}, \ldots,X_{u_{p}})|X_{s},X_{t}].$$
\end{lemma}

\noindent The proof of the previous lemma is elementary and will
be omitted.

\begin{lemma}
\label{lemmaA2} Let $(X, {\cal X}, \mu),(Y, {\cal Y}, \nu)$ be
probability spaces and $f:X \times Y \rightarrow {\bf R}$ a
bounded measurable function. If $\forall A \in {\cal X}, \forall B
\in {\cal Y}$
$$\int_{A \times B} f(x,y) (\mu \times \nu) (dx \times dy)= \mu(A) \int_{X \times B}f(x,y) (\mu \times \nu) (dx \times dy)$$
then $f(x,y)$ does not depend on $x$, for $(\mu \times
\nu)$-almost all $(x,y)$.
\end{lemma}

\noindent {\bf Proof}: Let $f_{B}(x)=\int_{B}f(x,y) \nu(dy), x \in
X$ and $I_{B}= \int_{X}f_{B}(x) \mu(dx)$. We have
$$\int_{A}f_{B}(x) \mu(dx)=\mu(A) I_{B} =\int_{A}I_{B} \mu(dx), \ \ \forall A \in {\cal X}$$
and hence $f_{B}(x)=I_{B}, \forall x \in N_{0}^{c}$, where $N_{0}$
is a $\mu$-negligible set. For each $x \in N_{0}^{c}$
$$\int_{B}f(x,y) \nu(dy)=\int_{X} \int_{B} f(x,y) \nu(dy) \mu(dx)= \int_{B} \int_{X}f(x,y) \mu(dx) \nu(dy), \ \forall B \in {\cal Y}.$$
Hence $f(x,y)=\int_{X}f(x,y) \mu(dx):=g(y)$ for all $y \in
N_{x}^{c}$, where $N_{x}$ is a $\nu$-negligible set. If we take
$N:=\{(x,y);x \in N_{0}^{c},y \in N_{x}^{c}\}^{c}$, then $(\mu
\times \nu)(N)=0$ and $f(x,y)=g(y), \forall (x,y) \in N^{c}$.
$\Box$

\section{A Bayes property of a Markov chain}

\begin{lemma}
\label{lemmaB1} Let $(Z_{1}, \ldots, Z_{k})$ be an increasing
Markov chain with values in $[0,1]$, with initial distribution
$\mu$ and transition probabilities $(Q_{i-1,i})_{i=2, \ldots,k}$;
let $Z_{0}:=0$ and $Z_{k+1}:=1$. Let $Y_{j}=Z_{j}-Z_{j-1};j=1,
\ldots,k+1$ and $X$ be a random variable such that
$${\cal P}[X=j|Y_{1}, \ldots,Y_{k+1}]=Y_{j} \ \ \ \forall j=1, \ldots,k+1.$$

\noindent Then for every $j=1, \ldots,k+1$, the conditional
distribution of $(Z_{1}, \ldots, Z_{k})$ given $X=j$ coincides
with the distribution of a Markov chain with some initial
distribution $\mu^{(j)}$ and some transition probabilities
$(Q_{i-1,i}^{(j)})_{i=2, \ldots,k}$.
\end{lemma}

\noindent {\bf Proof}: Let $\alpha_{j}:={\cal P}(X=j)={\cal
E}[Y_{j}]$. We consider first the case $j>1$. For any sets
$\Gamma_{1}, \ldots, \Gamma_{k} \in {\cal B}([0,1])$ we have
$${\cal P}[Z_{1} \in \Gamma_{1}, \ldots, \Gamma_{k} \in \Gamma_{k}|X=j] = \frac{1}{\alpha_{j}} \int_{\cap_{i=1}^{k}\{Z_{i} \in \Gamma_{i}\}} {\cal P}[X=j|Z_{1}, \ldots, Z_{k}]d {\cal P}=$$ $$\frac{1}{\alpha_{j}} \int_{\Gamma_{1}} \ldots \int_{\Gamma_{j}}h(z_{j})(z_{j}-z_{j-1})Q_{j-1,j}(z_{j-1};dz_{j}) \ldots Q_{12}(z_{1};dz_{2})\mu(dz_{1})$$

\noindent where $h(z_{j})= \int_{\Gamma_{j+1}} \ldots
\int_{\Gamma_{n}} Q_{k-1,k}(z_{k-1};dz_{k}) \ldots
Q_{j,j+1}(z_{j};dz_{j+1})$. We denote $\alpha_{j}^{(j)}(y) :={\cal
E}[Y_{j}|Z_{j-1}=y]$ and $\alpha_{i}^{(j)}(y):={\cal
E}[\alpha_{i+1}^{(j)}(Z_{i})|Z_{i-1}=y], i <j$; we have
$${\cal P}[Z_{i} \in \Gamma_{i}; i \leq k|X=j] =\frac{1}{\alpha_{j}} \int_{\Gamma_{1}} \alpha_{2}^{(j)}(z_{1}) \cdot \frac{1}{\alpha_{2}^{(j)}(z_{1})} \int_{\Gamma_{2}} \ldots \alpha_{j}^{(j)}(z_{j-1}) \cdot $$
$$\frac{1}{\alpha_{j}^{(j)}(z_{j-1})}
\int_{\Gamma_{j}}h(z_{j})(z_{j}-z_{j-1})Q_{j-1,j}(z_{j-1};dz_{j})
\ldots Q_{12}(z_{1};dz_{2}) \mu(dz_{1})=$$
$$\int_{\Gamma_{1}} \ldots \int_{\Gamma_{k}} Q_{k-1,k}^{(j)}(z_{k-1};dz_{k}) \ldots Q_{12}^{(j)}(z_{1};dz_{2}) \mu^{(j)}(dz_{1})$$
where $\mu^{(j)}(\Gamma):=(1/\alpha_{j})
\int_{\Gamma}\alpha_{2}^{(j)}(y) \mu(dy)$ and
\begin{eqnarray*}
Q_{i-1,i}^{(j)}(y; \Gamma) & := & Q_{i-1,i}(y; \Gamma) \ \ \ {\rm if} \ i>j \\
Q_{j-1,j}^{(j)}(y; \Gamma) & := & \frac{1}{\alpha_{j}^{(j)}(y)}\int_{\Gamma}(z-y)Q_{j-1,j}(y;dz) \\
Q_{i-1,i}^{(j)}(y; \Gamma) & := &
\frac{1}{\alpha_{i}^{(j)}(y)}\int_{\Gamma}
\alpha_{i+1}^{(j)}(z)Q_{i-1,i}(y;dz) \ \ \ {\rm if} \ i<j.
\end{eqnarray*}

\noindent We consider next the case $j=1$. For any sets
$\Gamma_{1}, \ldots, \Gamma_{k} \in {\cal B}([0,1])$ we have
$${\cal P}[Z_{1} \in \Gamma_{1}, \ldots, \Gamma_{k} \in \Gamma_{k}|X=1] = \frac{1}{\alpha_{1}} \int_{\cap_{i=1}^{k}\{Z_{i} \in \Gamma_{i}\}} {\cal P}[X=1|Z_{1}, \ldots,Z_{k}] d{\cal P}=$$ $$\int_{\Gamma_{1}} \ldots \int_{\Gamma_{k}} Q_{k-1,k}^{(1)}(z_{k-1};dz_{k}) \ldots Q_{12}^{(1)}(z_{1};dz_{2}) \mu^{(1)}(dz_{1})$$
where $\mu^{(1)}(\Gamma)=(1/\alpha_{1}) \int_{\Gamma}y \mu(dy)$
and $Q_{i-1,i}^{(1)}=Q_{i-1,i}, \forall i=2, \ldots,k$. $\Box$

\vspace{3mm}

The following definition is taken from Fang, Kotz and Ng (1990),
p.163.

\begin{definition}
\label{completely-neutral} {\rm A random vector $(Y_{1},
\ldots,Y_{k})$ with values in the simplex $S=\{(y_{j})_{j}; y_{j}
\in [0,1], \sum_{j=1}^{k}y_{j} \leq 1\}$ has a {\em completely
neutral} distribution if there exist some independent random
variables $V_{1}, \ldots,V_{k}$ such that $(Y_{1}, \ldots,Y_{k})$
has the same distribution as $(V_{1},V_{2}(1-V_{1}), \ldots, V_{k}
\prod_{j=1}^{k-1}(1-V_{j}))$.}
\end{definition}

The following result can be viewed as a complement to Theorem 4 of
Asgharian and Wolfson (2001).

\begin{corollary}
If $(Y_{1}, \ldots, Y_{k})$ has a completely neutral distribution,
$Y_{k+1}:=1- \sum_{j=1}^{k}Y_{j}$ and $X$ is a random variable
such that
$${\cal P}[X=j|Y_{1}, \ldots,Y_{k+1}]=Y_{j} \ \ \ \forall j=1, \ldots,k+1$$
then the conditional distribution of $(Y_{1}, \ldots, Y_{k})$
given $X=j$ is completely neutral.
\end{corollary}

\noindent {\bf Proof}: Let $Z_{i}=\sum_{j=1}^{i}Y_{j}$ and
$V_{1}=Y_{1}, V_{i}:=Y_{i}/(1-Z_{i-1}), i=2, \ldots,k$. The
variables $V_{1}, \ldots, V_{k}$ are independent and $Z_{1},
\ldots, Z_{k}$ is a Markov chain with the transition probabilities
$Q_{i-1,i}(y; \Gamma)=F_{i}((\Gamma-y)/(1-y))$, where $F_{i}$ is
the distribution of $V_{i}$. By Lemma \ref{lemmaB1}, the
conditional distribution of $(Z_{1}, \ldots, Z_{k})$ given $X=j$
coincide with the distribution of a Markov chain with some
transition probabilities $Q_{i-1,i}^{(j)}$. Direct calculations
show that $Q_{i-1,i}^{(j)}(y;
\Gamma)=F_{i}^{(j)}((\Gamma-y)/(1-y))$ with: $F_{i}^{(j)}=F_{i}$
if $i>j$,
$$F_{j}^{(j)}(\Gamma) = \frac{1}{\beta_{j}} \int_{\Gamma}v F_{j}(dv), \ \ F_{i}^{(j)}(\Gamma) = \frac{1}{1- \beta_{i}} \int_{\Gamma}(1-v) F_{i}(dv) \ \ \ {\rm if} \ i<j,$$
where $\beta_{i}={\cal E}[V_{i}]$. The conclusion follows
immediately. $\Box$

 }


\begin{thebibliography}{99}

\bibitem{asgharian-wolfson01} M. Asgharian and B.D. Wolfson, Covariates in multipath change-point problems: modelling and consistency of the maximum likelihood estimator, Canad. J. Statist. 29 (2001) 515--528.


\bibitem{balan-ivanoff02} R.M. Balan and B.G. Ivanoff, A Markov property for set-indexed processes, J. Theoret. Prob. 15 (2002) 553--588.

\bibitem{connor-mosimann69} R.J. Connor and J.E. Mosimann, Concepts of independence for proportions with a generalization of the Dirichlet distribution, J. Amer. Statist. Soc. 64 (1969) 194--206.

\bibitem{doksum74} K. Doksum, Tailfree and neutral random probabilities and their posterior distributions, Ann. Prob. 2 (1974) 183--201.

\bibitem{dykstra81} R.L. Dykstra and P. Laud, A Bayesian nonparametric approach to reliability, Ann. Statist. 9 (1981) 356--367.

\bibitem{fang-kotz-ng90} K. Fang, S. Kotz and K.W. Ng, Symmetric Multivariate and Related Distributions (1990), Chapman and Hall, London.

\bibitem{ferguson73} T.S. Ferguson, A Bayesian analysis of some nonparametric problems, Ann. Statist. 1 (1973) 209--230.

\bibitem{ferguson74} T.S. Ferguson, Prior distributions on spaces of probability measures, Ann. Statist. 2 (1974) 615--629.

\bibitem{ferguson-phadia74} T.S. Ferguson and E.G. Phadia, Bayesian nonparametric estimation based on censored data, Ann Statist. 7 (1979) 163--186.

\bibitem{halmos44} P.R. Halmos, Random alms, Ann. Math. Statist. 15 (1944) 182--189.

\bibitem{hjort90} N.L. Hjort, Nonparametric Bayes estimators based on Beta processes in models for life history data, Ann. Statist. 18 (1990) 1259--1294.

\bibitem{ivanoff-merzbach00} B.G. Ivanoff and E. Merzbach, Set-indexed Martingales (2000), Chapman and Hall, London.

\bibitem{walker-muliere97} S. Walker and P. Muliere, Beta-Stacy processes and a generalization of the P\'{o}lya-urn scheme, Ann. Statist. 25 (1997) 1762--1780.

\bibitem{walker-muliere99} S. Walker and P. Muliere, A characterization of a neutral to the right prior via an extension of Johnson's sufficientness postulate, Ann. Statist. 27 (1999) 589--599.

\end{thebibliography}
\end{document}